\RequirePackage{rotating}

\documentclass[12pt]{amsart}%
\usepackage{amsfonts}
\usepackage{amsmath}
\usepackage{amssymb}
\usepackage{amsthm}
\usepackage{booktabs}
\usepackage{mathtools}
\usepackage[colorlinks, linkcolor=red]{hyperref}
\usepackage[noabbrev,nameinlink]{cleveref}
\usepackage{enumerate}
\usepackage{bm}
\usepackage{tikz}
\usepackage{tikz-cd}
\usetikzlibrary{decorations.pathreplacing,patterns}
\usepackage{subcaption}
\usepackage{xcolor}
\usepackage{rotating}

\makeatletter
\newtheorem*{rep@theorem}{\rep@title}
\newcommand{\newreptheorem}[2]{%
\newenvironment{rep#1}[1]{%
 \def\rep@title{#2 \ref{##1}}%
 \begin{rep@theorem}}%
 {\end{rep@theorem}}}
\makeatother

\newtheorem{theorem}{Theorem}[section]
\newreptheorem{theorem}{Theorem}
\newtheorem{conjecture}[theorem]{Conjecture}
\newtheorem{corollary}[theorem]{Corollary}
\newtheorem{lemma}[theorem]{Lemma}
\newreptheorem{corollary}{Corollary}

\newtheorem{proposition}[theorem]{Proposition}

\theoremstyle{definition}
\newtheorem{definition}[theorem]{Definition}

\newtheorem{remark}[theorem]{Remark}
\newtheorem{example}[theorem]{Example}

\theoremstyle{plain}
\newtheorem*{namedthm}{\namedthmname}
\newcounter{namedthm}



\newcommand{\bC}{\mathbb{C}}

\newcommand{\bQ}{\mathbb{Q}}

\newcommand{\bZ}{\mathbb{Z}}
\newcommand{\cA}{\mathcal{A}}

\newcommand{\cD}{\mathcal{D}}
\newcommand{\cDH}{\mathcal{DH}}
\newcommand{\cDR}{\mathcal{DR}}

\newcommand{\cG}{\mathcal{G}}
\newcommand{\cH}{\mathcal{H}}
\newcommand{\cI}{\mathcal{I}}
\newcommand{\cJ}{\mathcal{J}}
\newcommand{\cK}{\mathcal{K}}

\newcommand{\cR}{\mathcal{R}}
\newcommand{\cS}{\mathcal{S}}

\newcommand{\cSH}{\mathcal{SH}}
\newcommand{\cSR}{\mathcal{SR}}

\newcommand{\fB}{\mathfrak{B}}
\newcommand{\fD}{\mathfrak{D}}
\newcommand{\fDih}{\mathfrak{Dih}}

\newcommand{\fS}{\mathfrak{S}}

\DeclareMathOperator{\Ann}{Ann}

\newcommand{\dif}{\mathrm{d}}
\DeclareMathOperator{\End}{End}
\DeclareMathOperator{\GL}{GL}

\DeclareMathOperator{\Hilb}{Hilb}
\DeclareMathOperator{\Hom}{Hom}

\DeclareMathOperator{\im}{im}
\DeclareMathOperator{\LT}{LT}

\DeclareMathOperator{\rS}{S}
\DeclareMathOperator{\sgn}{sgn}
\DeclareMathOperator{\Span}{Span}
\DeclareMathOperator{\Stir}{Stir}

\DeclareMathOperator{\vol}{vol}

\newcommand{\too}[1]{\overset{#1}{\to}}

\newenvironment{cd}[1][]
  {\begin{center}\begin{tikzcd}[#1]}
  {\end{tikzcd}\end{center}}
  
\newenvironment{textcolumn}
  {\begin{tabular}{@{}c@{}}}
  {\end{tabular}}

\title[Harmonic differential forms II. Bi-degree bounds]{Harmonic differential forms for pseudo-reflection groups \\ II. Bi-degree bounds}
\author{Joshua P. Swanson and Nolan R. Wallach}
\date{\today}

\begin{document}

\begin{abstract}
  This paper studies three results that describe the structure of the super-coinvariant algebra of pseudo-reflection groups over a field of characteristic $0$. Our most general result determines the top component in total degree, which we prove for all Shephard--Todd groups $G(m, p, n)$ with $m \neq p$ or $m=1$. Our strongest result gives tight bi-degree bounds and is proven for all $G(m, 1, n)$, which includes the Weyl groups of types $A$ and $B$/$C$. For symmetric groups (i.e. type $A$), this provides new evidence for a recent conjecture of Zabrocki \cite{1902.08966} related to the Delta Conjecture of Haglund--Remmel--Wilson \cite{MR3811519}. Finally, we examine analogues of a classic theorem of Steinberg \cite{MR167535} and the Operator Theorem of Haiman \cite{MR1918676,MR1256101}.

  Our arguments build on the type-independent classification of semi-invariant harmonic differential forms carried out in the first part of this series \cite{MR4258767}. In this paper we use concrete constructions including Gr\"obner and Artin bases for the classical coinvariant algebras of the pseudo-reflection groups $G(m, p, n)$, which we describe in detail. We also prove that exterior differentiation is exact on the super-coinvariant algebra of a general pseudo-reflection group. Finally, we discuss related conjectures and enumerative consequences.
\end{abstract}
\maketitle

\setcounter{tocdepth}{1}
\tableofcontents

\section{Introduction}

\subsection{Overview of results}
A polynomial $f \in \mathbb{Q}[x_1, \ldots, x_n]$ is $\fS_n$-\textit{harmonic} if it is annihilated by $\frac{\partial^k}{\partial x_1^k} + \cdots + \frac{\partial^k}{\partial x_n^k}$ for all $k \geq 1$. Here $\fS_n$ is the symmetric group on $n$ elements. A classic result of Steinberg \cite{MR167535} describes the harmonic polynomials as precisely those of the form $\partial_g \Delta_n$ where $\Delta_n \coloneqq \prod_{1 \leq i < j \leq n} (x_j - x_i)$ is the Vandermonde determinant and $\partial_g$ is a polynomial in the partial derivatives $\frac{\partial}{\partial x_i}$. In particular, the top-degree harmonic polynomial is $\Delta_n$, which transforms by the $\sgn$ representation of $\fS_n$.

Steinberg's result extends uniformly to an arbitrary \textit{pseudo-reflection group} $G$. We more generally consider the problem of determining the \textit{harmonic differential forms} of a pseudo-reflection group $G$. We have been motivated by Steinberg's result, Haiman's Operator Theorem \cite{MR1918676,MR1256101}, and a recent series of conjectures of Zabrocki \cite{1902.08966} and Haglund--Remmel--Wilson \cite{MR3811519} related to higher coinvariant algebras described in more detail below.

The first part of this series \cite{MR4258767} provides a starting point for such a description by giving a complete, type-independent construction of the $\det$-isotypic harmonic differential forms. In analogy with Steinberg's result, one may expect the $\det$-isotypic forms to be the ``top'' harmonics in some precise sense. Our main results are as follows:

\begin{enumerate}
  \item The top total-degree forms are $\det$-isotypic for ``almost all'' irreducible pseudo-reflection groups $G=G(m, p, n)$, namely those with $p \neq m$ or $m=1$ (\Cref{thm:C}).
  \item The top bi-degree forms are $\det$-isotypic for the groups $G(m, 1, n)$ (\Cref{thm:B}), which includes $\fS_n$ and the Coxeter groups of type $B$.
  \item All harmonic forms are obtained by applying partial derivatives to $\det$-isotypic forms when the rank is $\leq 2$ and either $G = G(m, 1, n)$ or $G$ is real (\Cref{thm:A1}) and also for multiples of the volume form when $G = G(m, 1, n)$ (\Cref{thm:A2}).
  \item The $t=0, z=-q$ specialization of the Hilbert series of the super-diagonal harmonics of the symmetric group agrees with Zabrocki's conjecture (see \Cref{ssec:intro_enumerative}).
\end{enumerate}

We also show that \Cref{thm:A1}, \Cref{thm:A2}, and \Cref{thm:B} fail for certain $G$; see \Cref{rem:A_fails} and \Cref{rem:B_fails}. We furthermore provide conjectures describing when they hold more broadly; see \Cref{conj:A} and \Cref{conj:C}.

Our arguments rely on the following tools which may be of independent interest.

\begin{enumerate}
  \setcounter{enumi}{4}
  \item The exterior derivative cochain complex on super-harmonic differential forms (or super coinvariants) is exact for all pseudo-reflection groups $G$ (\Cref{thm:exact}).
  \item The Artin and Gr\"obner bases for the coinvariant ideal for $G(m, p, n)$ (\Cref{ssec:intro_artin_grobner}, \Cref{thm:Artin_mpn}, \Cref{thm:grobner_mpn}).
\end{enumerate}

In the following subsections, we summarize classical properties of coinvariant algebras and harmonics for pseudo-reflection groups, introduce super analogues of these constructions, and state our main results and conjectures. See \Cref{sec:background} and \cite[\S2]{MR4258767} for additional detailed background.

\subsection{Coinvariant algebras and harmonics in type $A$}\label{ssec:intro_harmonics}

Let $\fS_n$ be the symmetric group on $n$ elements. The classical \textit{coinvariant algebra} of $\fS_n$ is the quotient
  \[ \cR_n \coloneqq \mathbb{Q}[x_1, \ldots, x_n]/\cI_n, \]
where $\cI_n \coloneqq \langle\mathbb{Q}[x_1, \ldots, x_n]_+^{\fS_n}\rangle$ is the \textit{coinvariant ideal} generated by all homogeneous symmetric polynomials of positive degree. A great deal is known about the structure of $\cR_n$ as a graded $\fS_n$-module \cite{MR526968}. The top-degree component of $\cR_n$ is the full subspace of elements which transform by $\sigma \cdot x = \sgn(\sigma)x$ and is spanned by the image of the Vandermonde determinant
  \[ \Delta_n \coloneqq \prod_{1 \leq i < j \leq n} (x_j - x_i). \]

The \textit{harmonics} of $\cI_n$ are
  \[ \cH_n \coloneqq \{f \in \mathbb{Q}[x_1, \ldots, x_n] : \partial_g f = 0 \text{ for all }g \in \cI_n\}. \]
Here $\partial_g$ is the partial differential operator defined by extending $x_i \mapsto \partial_{x_i} = \frac{\partial}{\partial x_i}$ multiplicatively and $\mathbb{Q}$-linearly. The harmonics $\cH_n$ are the orthogonal complement of $\cI_n$ under the natural positive-definite Hermitian form
\begin{equation}\label{eq:inner_product}
  (f, g) \coloneqq \text{degree zero component of }\partial_g f,
\end{equation}
so the natural projection $\cH_n \to \cR_n$ is a graded $\fS_n$-module isomorphism.

\subsection{Coinvariant algebras and harmonics of pseudo-reflection groups}\label{ssec:intro:psgs}

More generally, suppose $K \subset \mathbb{C}$ is a subfield closed under complex conjugation, $V = K^n$, and $G \subset \GL(V)$ is a \textit{pseudo-reflection group}. That is, $G$ is a finite group of unitary matrices generated by \textit{pseudo-reflections}, which are non-identity transformations of finite order which fix a hyperplane pointwise.

The pseudo-reflection groups were famously classified by Shephard--Todd \cite{MR0059914} into an infinite family $G(m, p, n)$ where $m, p, n \in \mathbb{Z}_{\geq 1}$ and $p \mid m$ together with $34$ exceptional groups. The group $G(1, 1, n) = \fS_n$ consists of $n \times n$ permutation matrices. The group $G(m, 1, n)$ consists of $n \times n$ \textit{pseudo-permutation matrices}, namely matrices such that each row and column has one non-zero entry which is an $m$th complex root of unity. In particular, $G(2, 1, n) = \fB_n$ is the Weyl group of signed permutations. Finally, the group $G(m, p, n)$ is the index-$p$ normal subgroup of $G(m, 1, n)$ consisting of pseudo-permutation matrices where the product of the non-zero elements taken to the $\frac{m}{p}$th power is $1$. In particular, $G(2, 2, n) = \fD_n$ is the Weyl group of signed permutations with an even number of signs.

The constructions in \Cref{ssec:intro_harmonics} generalize from $\fS_n$ to $G$ for any pseudo-reflection group $G$ (see \Cref{sec:background}). The \textit{coinvariant algebra of $G$} and the \textit{$G$-harmonics} are given as follows:
\begin{align*}
  \cR_G &\coloneqq K[x_1, \ldots, x_n]/\cI_G\qquad\text{where} \\
  \cI_G &\coloneqq \langle K[x_1, \ldots, x_n]_+^G\rangle\qquad\text{and} \\
  \cH_G &\coloneqq \{f \in K[x_1, \ldots, x_n] : \partial_g f = 0 \text{ for all }g \in \cI_G\},
\end{align*}
where $\partial_g$ is defined by extending $x \mapsto \partial_{x_i} = \frac{\partial}{\partial x_i}$ multiplicatively and \textit{conjugate}-linearly in $K$. For $G = G(m, p, n)$, the $G$-action is given by taking $x_i$ to be the $i$th coordinate function on $K^n$. As before, $\cH_G \to \cR_G$ is an isomorphism of graded $G$-modules.

Chevalley \cite{MR72877} showed that $\cR_G$, and hence $\cH_G$, carries the regular representation. Furthermore, the top-degree component of the harmonics $\cH_G$ is spanned by an element
  \[ \Delta_G \in K[x_1, \ldots, x_n] \]
which is unique up to non-zero scalar multiples and which transforms according to $g \cdot \Delta_G = \det(g) \Delta_G$. We call $\Delta_G$ the \textit{Vandermondian} of $G$, since $\Delta_{\fS_n} = \Delta_n$.

It is clear from the definition that if $f$ is harmonic, then so is $\partial_g f$.
Steinberg showed that every $G$-harmonic can be obtained ``from the top down'' starting with $\Delta_G$ as follows.

\begin{theorem}[Steinberg, {\cite[Thm.~1.3(c)]{MR167535}}]\label{thm:Steinberg}
  For any pseudo-reflection group $G$,
    \[ \cH_G = K[\partial_{x_1}, \ldots, \partial_{x_n}] \Delta_G = \{\partial_g \Delta_G : g \in K[x_1, \ldots, x_n]\}. \]
\end{theorem}

\subsection{Haiman's Operator Theorem}

A famous extension of the classical coinvariant algebra for $G = \fS_n$ was introduced by Garsia and Haiman \cite{MR1214091}. It involves two sets of commuting variables $x_1, \ldots, x_n$ and $y_1, \ldots, y_n$ with diagonal $\fS_n$-action given by $\sigma(x_i) \coloneqq x_{\sigma(i)}, \sigma(y_i) \coloneqq y_{\sigma(i)}$. The \textit{diagonal coinvariants} and \textit{diagonal harmonics} are $\fS_n$-modules bi-graded by $x$- and $y$-degree and are defined by
\begin{align*}
  \cDR_n &\coloneqq \mathbb{Q}[x_1, \ldots, x_n, y_1, \ldots, y_n]/\cK_n \qquad\text{where} \\
  \cK_n &\coloneqq \langle \mathbb{Q}[x_1, \ldots, x_n, y_1, \ldots, y_n]_+^{\fS_n}\rangle \qquad\text{and} \\
  \cDH_n &\coloneqq \{f \in \mathbb{Q}[x_1, \ldots, x_n, y_1, \ldots, y_n] : \partial_g f = 0 \text{ for all }g \in \cK_n\}.
\end{align*}

As usual, $\cDH_n \cong \cDR_n$ as bi-graded $\fS_n$-modules. Haiman conjectured and later proved the following description of the diagonal harmonics. Let $E_p \coloneqq \sum_{i=1}^n y_i \partial_{x_i}^p$. %One may check that the $E_p$ pairwise commute, and they commute with $\partial_{x_i}$. Indeed, they also preserve $\cDH_n$.

\begin{theorem}[{\cite[Thm.~4.2]{MR1918676}}]\label{thm:operator_theorem}
  We have
  \begin{equation}\label{eq:operator_theorem}
    \cDH_n = \mathbb{Q}[\partial_{x_1}, \ldots, \partial_{x_n}, E_1, \ldots, E_{n-1}] \Delta_n.
  \end{equation}
\end{theorem}

\noindent Haiman's proof of this theorem, which was originally conjectured in \cite[Conj.~5.1.1]{MR1256101}, involves the deep use of algebraic geometry. An elementary proof of the $y$-degree $1$ component of \Cref{thm:operator_theorem} was given by Alfano \cite{MR1661359}.

\subsection{Super coinvariant algebras and harmonic differential forms}\label{ssec:intro_super_harmonics}

A recent conjecture of Zabrocki \cite{1902.08966} introduced the \textit{super-diagonal coinvariant algebra} as a potential representation-theoretic model for the \textit{Delta Conjecture} of Haglund--Remmel--Wilson \cite{MR3811519}. The $t=0$ case of Zabrocki's conjecture sparked significant interest in the following extension of the classical coinvariant algebras, which is our main object of interest. See \cite[\S1.2]{MR4258767} for an overview of this and related work and further references. This paper studies the extension of the polynomial ring $K[x_1, \ldots, x_n]$ by adjoining \textit{anti-commuting} variables $\theta_1, \ldots, \theta_n$ where $\theta_i \theta_j = -\theta_j \theta_i$ and $x_i \theta_j = \theta_j x_i$. Here $\theta_i$ is conceptually the differential $1$-form $\dif x_i$, so the products $\theta_{i_1} \cdots \theta_{i_k}$ represent differential $k$-forms $\dif x_{i_1} \wedge \cdots \wedge \dif x_{i_k}$. We use the $\theta$ variables for consistency with existing literature. Since $\theta_i^2 = 0$, we often take $i_1 < \cdots < i_k$. Let $K[x_1, \ldots, x_n, \theta_1, \ldots, \theta_n]$ be the ring of differential forms with polynomial coefficients. The $G$-action on the $\theta_i$ is the same as on the $x_i$ and is extended multiplicatively and $K$-linearly to $K[x_1, \ldots, x_n, \theta_1, \ldots, \theta_n]$, which is a $G$-module bi-graded by $x$-degree and $\theta$-degree. Abstractly, this is the ring $\rS(V^*) \otimes \wedge V^*$; see \cite[\S2]{MR4258767} for details.

\begin{definition}
  The \textit{super coinvariant algebra} of a pseudo-reflection group $G$ is the quotient
  \[ \cSR_G \coloneqq K[x_1, \ldots, x_n, \theta_1, \ldots, \theta_n]/\cJ_G \]
  where $\cJ_G \coloneqq \langle K[x_1, \ldots, x_n, \theta_1, \ldots, \theta_n]_+^G\rangle$ is the \textit{super coinvariant ideal} generated by bi-homogeneous $G$-invariant differential forms of positive total degree.
\end{definition}

We write $\cSR_G^k$ for the component of $\cSR_G$ consisting of images of $k$-forms, i.e. the component of $\theta$-degree $k$. Note that $\cSR_G^0 = \cR_G$, so the super coinvariant algebra of $G$ contains the classical coinvariant algebra of $G$.

To define the harmonics in this context requires the extension of the partial differential operators $\partial_g$ to differential forms. Let $\partial_{\theta_i}$ be the usual \textit{interior product} defined by 
\begin{align*}
  \partial_{\theta_i} \theta_{i_1} \cdots \theta_{i_k}
    &= \begin{cases}
            (-1)^{\ell-1} \theta_{i_1} \cdots \widehat{\theta_{i_\ell}} \cdots \theta_{i_k}
              & \text{if }i = i_\ell \\
            0 & \text{otherwise}.
          \end{cases}
\end{align*}
Here $\partial_{x_i}$ is $\theta$-linear and $\partial_{\theta_i}$ is $x$-linear. Note that $\partial_{x_i}$ and $\partial_{\theta_j}$ satisfy the same (anti-)commutation relations as the $x_i$ and $\theta_j$. For $\omega \in K[x_1, \ldots, x_n, \theta_1, \ldots, \theta_n]$, let $\partial_\omega$ be obtained by replacing each $x_i$ with $\partial_{x_i}$, replacing each $\theta_j$ with $\partial_{\theta_j}$, reversing the order of the $\theta$'s, and taking the conjugate of the coefficients. These twists ensure the extension of \eqref{eq:inner_product} remains positive-definite.

\begin{definition}
  The \textit{super harmonics} of a pseudo-reflection group $G$ are
  \[ \cSH_G \coloneqq \{\eta \in K[x_1, \ldots, x_n, \theta_1, \ldots, \theta_n] : \partial_\omega \eta = 0\text{ for all }\omega \in \cJ_G\}. \]
\end{definition}
\noindent The natural projection $\cSH_G \to \cSR_G$ is again an isomorphism of bi-graded $G$-modules, so $\cSR_G$ and $\cSH_G$ are frequently interchangeable.

\subsection{The $\det$-isotypic harmonic differential forms}\label{ssec:intro_semi-invariants}

From Steinberg's \Cref{thm:Steinberg}, the $\det$-isotypic component of the harmonics,
  \[ \cH_G^{\det} \coloneqq \{f \in \cH_G : \sigma \cdot f = \det(\sigma)f\} = \Span_K\{\Delta_G\}, \]
plays a special role. In \cite{MR4258767}, the authors gave a ``top-down'' construction of $\cSH_G^{\det}$ in the spirit of Steinberg's \Cref{thm:Steinberg} and Haiman's \Cref{thm:operator_theorem} using certain differential operators $\dif_1, \ldots, \dif_r \in \End_K(\cSH_G)$. Here $r \coloneqq \dim(V/V^G)$. We call these operators \textit{generalized exterior derivatives} since in general we may take
\begin{equation}\label{eq:exterior_derivative}
  \dif_1 = \dif \coloneqq \sum_{i=1}^n \frac{\partial}{\partial x_i} \theta_i
\end{equation}
to be the \textit{exterior derivative}, where $\theta_i$ here denotes left multiplication by $\theta_i$. If $e_1^*, \ldots, e_r^*$ are the positive \textit{co-exponents} of $G$, then $\dif_i$ lowers $x$-degree by $e_i^*$ and raises $\theta$-degree by $1$. See \Cref{ssec:generalized_exterior_derivatives} for details and \Cref{tab:explicit} for explicit descriptions of $\dif_1, \ldots, \dif_r$ for $G(m, p, n)$. \Cref{tab:explicit} lists the real cases $G = \fS_n, \fB_n, \fD_n, \fDih_{2n}$ explicitly.

The main result of \cite{MR4258767} gives a basis of $2^r$ elements for the $\det$-isotypic component of the super harmonics.

\begin{theorem}[{\cite[Thm.~5.7]{MR4258767}}]\label{thm:semi-invariants}
  Let $G$ be a pseudo-reflection group. Then
  \begin{equation}\label{eq:SH_det_basis}
    \cSH_G^{\det} = \Span_K\{\dif_{i_1} \cdots \dif_{i_k} \Delta_G : 1 \leq i_1 < \cdots < i_k \leq r\}.
  \end{equation}
\end{theorem}

\subsection{Differential operator results}

The partial derivatives $\partial_g$ for $g \in K[x_1, \ldots, x_n]$ clearly preserve $\cSH_G$, so $\partial_g \cSH_G^{\det} \subset \cSH_G$. Motivated by Steinberg's \Cref{thm:Steinberg} and Haiman's \Cref{thm:operator_theorem}, we show that the following reverse containments hold.

\begin{theorem}\label{thm:A1}
  Let $G \subset \GL(V)$ be a pseudo-reflection group with rank $r = \dim(V/V^G)$. Then if $r \leq 2$ and either $G = G(m, 1, n)$ or $G$ is real,
  \begin{equation}\label{eq:thm:A}
  \begin{split}
    \cSH_G
      &= K[\partial_{x_1}, \ldots, \partial_{x_n}] \cSH_G^{\det} \\
      &= \{\partial_g \dif_{i_1} \cdots \dif_{i_k} \Delta_G : g \in K[x_1, \ldots, x_n], i_j \in [r] \}.
  \end{split}
  \end{equation}
\end{theorem}

\begin{theorem}\label{thm:A2}
  Let $G = G(m, 1, n)$ or let $G$ be real. Then the $\theta$-degree $r$ component of \eqref{eq:thm:A} holds.
\end{theorem}

\begin{remark}\label{rem:A_fails}
  In \Cref{lem:no_dice}, we show that in fact the $\cSH_G^r$ component of \eqref{eq:thm:A} holds for $G=G(m, p, n)$ if and only if $r \leq 2$, $G=G(m, 1, n)$, or $G = G(2, 2, n)$. Hence \eqref{eq:thm:A} cannot possibly hold for any groups in the infinite family $G(m, p, n)$ beyond these. However, computational data shows that \eqref{eq:thm:A} \textit{fails} for $G = \fD_4$ and $\fD_5$. It appears likely to fail for $\fD_n$ with $n \geq 4$. On the other hand, we have verified that \eqref{eq:thm:A} holds for $\fS_n$ with $n \leq 6$, $\fB_n$ with $n \leq 4$, $G(3, 1, 4)$, and $G(5, 1, 3)$, among others. See \Cref{tab:calcs} for additional data. We also note that the $\cSH_{\fS_n}^1$ case of \eqref{eq:thm:A} is equivalent to the $y$-degree $1$ case of Haiman's \Cref{thm:operator_theorem}. Consequently, we are lead to the following conjecture.
\end{remark}

\begin{conjecture}[Differential Operator Conjecture]\label{conj:A}
  If $G = G(m, 1, n)$, then \eqref{eq:thm:A} holds.
\end{conjecture}

\Cref{thm:A1} includes the dihedral groups $G(m, m, 2) = \fDih_{2m}$ ($m \geq 1$) and $6$ exceptional groups. We do not have a complete determination or conjecture for the exceptional groups for which \eqref{eq:thm:A} holds. It does hold for $H_3$, though perhaps surprisingly it fails for $F_4$.

Our proof of \Cref{thm:A1} is mostly uniform and relies on the following result which may be of independent interest. See \Cref{ssec:intro_enumerative} for additional consequences.

\begin{theorem}\label{thm:exact}
  For any pseudo-reflection group $G \subset \GL(V)$ with $r = \dim(V/V^G)$, the exterior derivative cochain complex
  \begin{equation}\label{eq:exact_d}
    0 \to K \to \cSH_G^0 \too{\dif} \cSH_G^1 \too{\dif} \cdots \too{\dif} \cSH_G^r \too{\dif} 0
  \end{equation}
  is exact.
\end{theorem}

\noindent The complex in \Cref{thm:exact} is a finite-dimensional, algebraic analogue of the de Rham complex of a smooth manifold. Exactness is proved with an analogue of Hodge theory using total Laplacians. See \Cref{sec:complex} for the proofs of \Cref{thm:A1} and \Cref{thm:exact}.

\subsection{Bi-degree bound results}

When \eqref{eq:thm:A} holds, the harmonics
  \[ \dif_{i_1} \cdots \dif_{i_k} \Delta_G \]
are the ``top-most'' elements of $\cSH_G$. In particular, it implies that for each $k$, the top $x$-degree elements of $\cSH_G$ belong to $\cSH_G^{\det}$ and fully describes the non-zero bi-degree components of $\cSH_G$ as follows. Write $\cSH_G^{i, k}$ for the $x$-degree $i$, $\theta$-degree $k$ component of $\cSH_G$.

\begin{lemma}\label{lem:B}
  Let $G \subset \GL(V)$ have rank $r = \dim(V/V^G)$ and co-exponents $1 \leq e_1^* \leq \cdots \leq e_r^*$. Suppose \eqref{eq:thm:A} holds. Then
  \begin{equation}\label{eq:lem:B}
  \begin{split}
    \cSH_G^{i, k} \neq 0 \quad\Leftrightarrow\quad &0 \leq k \leq r \text{ and } \\
        &0 \leq i \leq \deg \Delta_G - e_1^* - \cdots - e_k^*.
  \end{split}
  \end{equation}
  Moreover, if $i = \deg \Delta_G - e_i^* - \cdots - e_k^*$, then $\cSH_G^{i, k} \subset \cSH_G^{\det}$.
\end{lemma}

Our strongest result is to show that the following consequence of \Cref{conj:A} is true unconditionally. In particular, it verifies the predicted bi-degree support of Zabrocki's super coinvariant algebra conjecture when $t=0$, providing additional evidence for that conjecture.

\begin{theorem}\label{thm:B}
  Let $G = G(m, 1, n)$. Then $\cSH_{G(m, 1, n)}^{i, k} \neq 0$ if and only if $i, k \geq 0$ and
  \begin{equation}\label{eq:Gm1n_support}
    i+k+m\binom{k}{2} \leq m \binom{n}{2} + (m-1)n.
  \end{equation}
%  Moreover, if equality holds in \eqref{eq:Gm1n_support}, then
%    \[ \cSH_{G(m, 1, n)}^{i, k} = K \dif_1 \cdots \dif_r \prod_{1 \leq a < b \leq n} (x_b^m - x_a^m) \cdot (x_1 \cdots x_n)^{m-1}. \]
  
%  \begin{alignat*}{3}
%    \cSH_{G(m, 1, n)}^{i, k} \neq 0 \quad\Leftrightarrow\quad &0 \leq k \leq r \\
%       &\text{and } 0 \leq i &&\leq m\binom{n}{2} + n(m-1) - m\binom{k}{2} - k \\
%          &&&= \frac{(n+k)(m(n-k+1)-2)}{2}.
%  \end{alignat*}
%  Moreover, if $i = \frac{(n+k)(m(n-k+1)-2)}{2}$, then
%    \[ \cSH_{G(m, 1, n)}^{i, k} = K \dif_1 \cdots \dif_r \prod_{1 \leq i < j \leq n} (x_j^m - x_i^m) \cdot (x_1 \cdots x_n)^{m-1}. \]
\end{theorem}

%We have also verified the \ref{conj:B} when $G$ is the Weyl group of type $E_6$, $F_4$, or $D_n$ for $n \leq 8$, and the rank $2$ case follows from \Cref{thm:A1}. The remaining cases are hence the Weyl groups of type $H_3$, $H_4$, $E_7$, $E_8$, and $D_n$ for $n \geq 9$.

\begin{remark}\label{rem:B_fails}
  Note that \Cref{lem:B} may hold even when \eqref{eq:thm:A} fails. By \Cref{tab:calcs}, this indeed occurs for $\fD_4, \fD_5, F_4$. In \Cref{lem:no_dice}, we show that the $\cSH_G^r$ component of \eqref{eq:lem:B} holds for $G(m, p, n)$ if and only if $r \leq 2$, $G = G(m, 1, n)$, or $G = G(2, 2, n) = \fD_n$. In contrast to \Cref{conj:A}, our data do not rule out the possibility that \eqref{eq:lem:B} holds for $\fD_n$.
\end{remark}

Our proof of \Cref{thm:B} uses Gr\"{o}bner and Artin bases of $\cR_{G(m, 1, n)}$ developed in \Cref{sec:artin}. See \Cref{sec:bidegree}.

\subsection{Total degree bound results}

The bi-degree support and top components from \Cref{lem:B} imply the following total degree support and top components of $\cSH_G$, or equivalently $\cSR_G$.

\begin{lemma}\label{lem:C}
  Let $G$ be a pseudo-reflection group. Suppose \Cref{lem:B} holds. Then
  \begin{equation}\label{eq:lem:C.1}
    \bigoplus_{i+k = d} \cSH_G^{i, k} \neq 0 \quad\Leftrightarrow\quad
      0 \leq d \leq \deg \Delta_G.
  \end{equation}
  Moreover,
  \begin{equation}\label{eq:lem:C.2}
    \bigoplus_{i+k=\deg \Delta_V} \cSH_G^{i, k} \subset \cSH_G^{\det}.
  \end{equation}
\end{lemma}

We show that this weaker description is true even in many cases where \Cref{lem:B} fails to hold. Our most general result is the following.

\begin{theorem}\label{thm:C}
  Let $G = G(m, p, n)$ with $p \neq m$ or $p=1$. Then
    \[ \bigoplus_{i+k=\ell} \cSH_{G(m, p, n)}^{i, k} \neq 0 \quad\Leftrightarrow\quad
        0 \leq \ell \leq m \binom{n}{2} + n \left(\frac{m}{p} - 1\right). \]
  Moreover, if $\ell = m \binom{n}{2} + n \left(\frac{m}{p} - 1\right)$, then
    \[ \bigoplus_{i+k=\ell} \cSH_{G(m, p, n)}^{i, k} = 
    (K + K\dif) \prod_{1 \leq i < j \leq n} (x_j^m - x_i^m) \cdot (x_1 \cdots x_n)^{m/p - 1}. \]\end{theorem}

Our proof of \Cref{thm:C} again uses Gr\"{o}bner and Artin bases of $\cR_{G(m, 1, n)}$ developed in \Cref{sec:artin}. See \Cref{sec:total}.

Our results and computations have uncovered no cases in which \Cref{thm:C} fails to hold. Indeed, our argument shows that it must hold for all $G = G(m, p, n)$ except possibly when $m=p>1$ and $k \in \{n-1, n-2\}$. Consequently, we conjecture the following.

\begin{conjecture}[Total Degree Bounds Conjecture]\label{conj:C}
  Let $G$ be a pseudo-reflection group. Then
    \[ \bigoplus_{i+k = d} \cSH_G^{i, k} \neq 0 \quad\Leftrightarrow\quad
        0 \leq d \leq \deg \Delta_G. \]
  Moreover, $\bigoplus_{i+k=\deg \Delta_V} \cSH_G^{i, k} \subset \cSH_G^{\det}$.
\end{conjecture}

\subsection{Hilbert series considerations}\label{ssec:intro_enumerative}

We finish this introduction with some additional enumerative consequences of the preceding results which support some further conjectures.

The \textit{Hilbert series} of $\cSR_G^k$ and $\cSR_G$ are
\begin{align*}
  \Hilb(\cSR_G^k; q) &\coloneqq \sum_{i \geq 0} q^i \dim_K(\cSR_G^{i, k}) \\
  \Hilb(\cSR_G; q, z) &\coloneqq \sum_{k \geq 0} z^k \Hilb(\cSR_G^k; q),
\end{align*}
where $q$ tracks $x$-degree and $z$ tracks $\theta$-degree. An immediate consequence of \Cref{thm:exact} is the following enumerative corollary.

\begin{corollary}\label{cor:Hilb}
  For any pseudo-reflection group $G$,
    \[ \Hilb(\cSR_G; q, -q) = \sum_{k \geq 0} (-q)^k \Hilb(\cSR_G^k; q) = 1. \]
\end{corollary}

Our overarching goal has been to provide evidence for Zabrocki's conjecture \cite{1902.08966} for the tri-graded Frobenius series of the type $A$ super-diagonal coinvariant algebra. In particular, we may check Zabrocki's conjecture against \Cref{cor:Hilb}.

To do so, let $\Stir_q(n, k)$ be a \textit{$q$-Stirling number of the second kind} \cite[\S3]{MR27288}, defined recursively by $\Stir_q(n, k) = \Stir_q(n-1, k-1) + [k]_q \Stir_q(n-1, k)$ and $\Stir_q(1, k) = \delta_{k=1}$. Here $[k]_q \coloneqq 1 + q + \cdots + q^{k-1}$ is a \textit{$q$-integer} and $\delta_P = 1$ if $P$ is true and $0$ otherwise. Zabrocki's conjecture specializes as follows.

\begin{conjecture}[Zabrocki, {\cite{1902.08966}}]\label{conj:Hilb_type_A}
  For $0 \leq k \leq n-1$,
    \[ \Hilb(\cSR_{\fS_n}^k; q) = [n-k]_q! \Stir_q(n, n-k). \]
\end{conjecture}

\noindent One may indeed check that
  \[ \sum_{k=0}^{n-1} (-q)^k [n-k]_q! \Stir_q(n, n-k) = 1, \]
so \Cref{conj:Hilb_type_A} is consistent with \Cref{cor:Hilb}.

Based on computational evidence including \Cref{tab:calcs}, the first named author has introduced a type $B$ analogue of \Cref{conj:Hilb_type_A}. Let $\Stir_q^B(n, k)$ be a \textit{type $B$ $q$-Stirling number of the second kind}, defined recursively by $\Stir_q^B(n, k) = \Stir_q^B(n-1, k-1) + [2k+1]_q \Stir_q^B(n-1, k)$ and $\Stir_q^B(0, k) = \delta_{k=0}$. The type $B$ analogue of \Cref{conj:Hilb_type_A} is as follows.

\begin{conjecture}\label{conj:Hilb_type_B}
  For $0 \leq k \leq n$,
    \[ \Hilb(\cSR_{\fB_n}^k; q) = [n-k]_{q^2}! [2]_q^{n-k} \Stir_q^B(n, n-k). \]
\end{conjecture}

\noindent One may again check that
  \[ \sum_{k=0}^n (-q)^k [n-k]_{q^2}! [2]_q^{n-k} \Stir_q^B(n, n-k) = 1, \]
so \Cref{conj:Hilb_type_B} is also consistent with \Cref{cor:Hilb}. An investigation of the combinatorial properties of these and other $q$-analogues is in progress \cite{sss21}.

One may obtain a different complex from \eqref{eq:exact_d} by replacing the differentials $\dif$ with $\dif_i$ for $1 \leq i \leq r$, though the result is typically not exact. The graded Euler characteristic of the complex is
  \[ \chi(H^*(\cSH_G, \dif_i); q) = \Hilb(\cSR_G; q, -q^{e_i^*}) - 1. \]
A potential approach to \Cref{conj:Hilb_type_A} and \Cref{conj:Hilb_type_B} is to find homotopic complexes with the correct Euler characteristic. More concretely, \Cref{conj:Hilb_type_A} is equivalent to the following variation on \Cref{cor:Hilb}.

\begin{conjecture}\label{conj:Hilb_alt}
  For $1 \leq j \leq n-1$,
    \[ \Hilb(\cSR_{\fS_n}; q, -q^j) = \sum_{k=0}^n (-q^j)^k [n-k]_q! \Stir_q(n, n-k). \]
\end{conjecture}

%In pursuit of further evidence for Zabrocki's conjecture, we consider generalizing \Cref{thm:B} to the super diagonal coinvariants in \Cref{sec:tridegree}. We reduce the relevant bound to a conjecture involving two sets of commuting variables, \Cref{conj:twocomm}, and provide computational evidence for this conjecture.

\subsection{Paper organization}

%The second author previously presented the heart of the argument for \Cref{thm:B} in \cite{1906.11787}. The present work extends, generalizes, and supercedes \cite{1906.11787}.

\Cref{sec:background} gives background on polynomial differential forms, generalized exterior derivatives, and the invariant theory of pseudo-reflection groups. \Cref{sec:top} analyzes the structure of the top $\theta$-degree component of $\cSH_G$ and proves \Cref{thm:A2}. \Cref{sec:complex} set up our algebraic Hodge theory argument proving exactness of exterior differentiation on the harmonics, \Cref{thm:exact}, as well as \Cref{thm:A1}. In \Cref{sec:artin} we describe the Artin and Gr\"obner bases for $G(m, p, n)$; see \Cref{thm:Artin_mpn} and \Cref{thm:grobner_mpn}. \Cref{sec:bidegree} proves the bi-degree bounds in \Cref{thm:B}. \Cref{sec:total} considers the top total degree and proves \Cref{thm:C}.
%\Cref{sec:tridegree} considers tri-degree bounds in Zabrocki's conjecture.

\section{Background}\label{sec:background}

\subsection{Polynomial differentials}

We first briefly introduce the standard $G$-module structures and differential operators underlying our results. See \cite[\S2]{MR4258767} for a general, abstract version. The following concrete version is included in the spirit of much of the combinatorics literature and is intended to make these developments more accessible.

Let $K \subset \mathbb{C}$ be a subfield closed under complex conjugation. Let $V = K^n$. Suppose $G \subset U(n, K)$ is a subgroup of unitary matrices, so $(\sigma^{-1})_{ij} = \overline{\sigma}_{ji}$ with $\sigma_{ij}$ defined by $\sigma \cdot f_j = \sum_{i=1}^n \sigma_{ij} f_i$ where $f_1, \ldots, f_n$ is the standard orthonormal basis. Let $x_1, \ldots, x_n \in V^* \coloneqq \Hom_K(V, K)$ be the dual basis $x_i(f_j) = \delta_{i=j}$.

The ring $K[x_1, \ldots, x_n]$ consists of the polynomial functions $f \colon V \to K$. The group $G$ acts naturally on polynomial functions via the \textit{contragredient} action $(\sigma \cdot f)(v) \coloneqq f(\sigma^{-1}(v))$. Concretely, $\sigma \cdot x_j = \sum_{i=1}^n \overline{\sigma}_{ij} x_i$ and $\sigma \cdot \sum c_\alpha x^\alpha = \sum c_\alpha \sigma(x_1)^{\alpha_1} \cdots \sigma(x_n)^{\alpha_n}$. Define a \textit{conjugate}-linear bijection
\begin{align*}
  \tau \colon V \to V^* \qquad\text{where}\qquad
  \sum_{i=1}^n c_i f_i \mapsto \sum_{i=1}^n \overline{c}_i x_i
\end{align*}
Then $\tau$ is $G$-equivariant, i.e.~$\tau(\sigma \cdot v) = \sigma \cdot \tau(v)$ for all $\sigma \in G$.

The derivative of $f \in K[x_1, \ldots, x_n]$ in the direction $v \in V$ is the polynomial function $\partial_v f$ defined by $(\partial_v f)(w) \coloneqq \left.\frac{\dif}{\dif t}\right|_{t=0} f(w+tv)$. To simplify the exposition, we transfer these derivatives to $V^*$ by defining operators $\partial_x f \coloneqq \partial_{\tau^{-1}(x)} f$. Concretely, $\partial_{x_i} = \frac{\partial}{\partial x_i}$ is the usual partial derivative. More generally, we extend these partial derivatives to polynomials in $K[x_1, \ldots, x_n]$ multiplicatively and \textit{conjugate}-linearly.

\begin{definition}
  If $g = \sum_\alpha c_\alpha x^\alpha \in K[x_1, \ldots, x_n]$ then set
  \[ \partial_g \coloneqq \sum_\alpha \overline{c}_\alpha \partial_{x_1}^{\alpha_1} \cdots \partial_{x_n}^{\alpha_n}. \]
\end{definition}
We use two fundamental properties of $\partial_g$. First,
$\sigma \cdot \partial_g f = \partial_{\sigma \cdot g}(\sigma \cdot f)$. In particular, if $g$ is $G$-invariant, then $\sigma \cdot \partial_g f = \partial_g(\sigma \cdot f)$ and $\partial_g$ is $G$-equivariant. Second, if $g$ is homogeneous and non-zero, then $\partial_g g = \sum_\alpha |c_\alpha|^2 \alpha! > 0$.

We may define a $G$-invariant positive-definite Hermitian form on $K[x_1, \ldots, x_n]$ by
  \[ (f, g) \coloneqq (\partial_g f)(0, \ldots, 0). \]
This form is linear in $f$ and conjugate-linear in $g$. Under this form, the coinvariant ideals and harmonics from \Cref{ssec:intro_harmonics} are orthogonal complements, $\cI_G^\perp = \cH_G$ and $\cH_G^\perp = \cI_G$. Furthermore, the adjoint of $\partial_{x_i}$ with respect to this form is multiplication by $x_i$, i.e.~$\partial_{x_i}^\dagger f = x_i f$. More generally, $\partial_g^\dagger f = gf$. See for example \cite[Lemmas~2.7, 2.11, 4.2]{MR4258767}.

\subsection{Differential forms}\label{ssec:differential_forms}

In place of polynomial functions $f(v)$ on $V$, we may consider alternating multilinear $k$-forms $\eta \colon V^k \to K$ on $V$. We again have the natural contragredient action $(\sigma \cdot \eta)(v_1, \ldots, v_k) \coloneqq \eta(\sigma^{-1}(v_1), \ldots, \sigma^{-1}(v_k))$. Concretely, the alternating $k$-forms on $V$ have $\binom{n}{k}$ basis elements indexed by $1 \leq i_1 < \cdots < i_k \leq n$ and determined by
  \[ (\dif x_{i_1} \wedge \cdots \wedge \dif x_{i_k})(f_{j_1}, \ldots, f_{j_k})
      = \begin{cases}
        1 & \text{if }j_1 = i_1, \ldots, j_k = i_k \\
        0 & \text{if }\{i_1, \ldots, i_k\} \neq \{j_1, \ldots, j_k\}.
      \end{cases} \]

Following a standard convention in this area of algebraic combinatorics, we abbreviate $\theta_i \coloneqq \dif x_i$ and typically suppress the $\wedge$ symbol. The alternating multilinear forms on $V$ under the wedge product together form the $K$-algebra $K[\theta_1, \ldots, \theta_n]$ generated by $\theta_1, \ldots, \theta_n$ subject to the relations $\theta_i\theta_j = -\theta_j\theta_i$, and in particular $\theta_i^2 = 0$. The $G$-action is given by $\sigma \cdot \theta_j = \sum_{i=1}^n \overline{\sigma}_{ij} \theta_i$ and $\sigma \cdot \sum c_{i_1, \ldots, i_k} \theta_{i_1} \cdots \theta_{i_k} = \sum c_{i_1, \ldots, i_k} \sigma(\theta_{i_1}) \cdots \sigma(\theta_{i_k})$.

The analogue of the directional derivative $\partial_v$ for alternating forms is given by the \textit{interior product}. If $\eta \colon V^k \to K$ is alternating and $v \in V$, then $(\partial_v \eta)(v_1, \ldots, v_{k-1}) \coloneqq \eta(v, v_1, \ldots, v_{k-1})$. In contrast to partial derivatives, which commute, interior products anti-commute according to $\partial_v \partial_w = -\partial_w \partial_v$. We again have $\sigma \cdot \partial_v\eta = \partial_{\sigma(v)}(\sigma \cdot \eta)$. By a slight abuse of notation, now let $\tau \colon V \to V^*$ be the conjugate-linear bijection given by $\sum_{i=1}^n c_i f_i \mapsto \sum_{i=1}^n \overline{c}_i \theta_i$, which is again $G$-equivariant. Let $partial_\theta \coloneqq \partial_{\tau^{-1}(\theta)}$. Concretely, if $i_1 < \cdots < i_k$,
  \[ \partial_{\theta_i} \theta_{i_1} \cdots \theta_{i_k}
      = \begin{cases}
      (-1)^{\ell-1} \theta_{i_1} \cdots \widehat{\theta}_{i_\ell} \cdots \theta_{i_k} & \text{if }i = i_\ell \\
      0 & \text{otherwise},
      \end{cases} \]
where $\widehat{\theta}_{i_\ell}$ means $\theta_{i_\ell}$ is omitted. 

\begin{definition}
Let $K[x_1, \ldots, x_n, \theta_1, \ldots, \theta_n]$ be the ring of differential forms on $V$ with coefficients which are polynomial functions on $V$. This is the $K$-algebra generated by $x_1, \ldots, x_n, \theta_1, \ldots, \theta_n$ subject to the relations $x_i x_j = x_j x_i$, $\theta_i \theta_j = -\theta_j \theta_i$, and $x_i \theta_j = \theta_j x_i$.
\end{definition}

As before, we extend the interior product operators multiplicatively and conjugate-linearly to $K[x_1, \ldots, x_n, \theta_1, \ldots, \theta_n]$. Here we write $I$ in place of $1 \leq i_1 < \cdots < i_k \leq n$. Note the reversal of the order of the $\partial_\theta$ operators.

\begin{definition}
  If
    \[ \omega = \sum_{\alpha, I} c_{\alpha, I} x^\alpha \theta_{i_1} \cdots \theta_{i_k} \in K[x_1, \ldots, x_n, \theta_1, \ldots, \theta_n] \]
  then set
  \begin{equation}\label{eq:partial_omega}
    \partial_\omega \coloneqq \sum_{\alpha, I} \overline{c}_{\alpha, I} \partial_{x_1}^{\alpha_1} \cdots \partial_{x_n}^{\alpha_n} \partial_{\theta_{i_k}} \cdots \partial_{\theta_{i_i}}.
  \end{equation}
\end{definition}

Note that $\sigma \cdot \partial_\omega f = \partial_{\sigma \cdot \omega} (\sigma \cdot f)$. Moreover, one may check that if $\omega$ is bi-homogeneous and non-zero, then $\partial_\omega \omega = \sum_{\alpha, I} |c_{\alpha, I}|^2 \alpha! > 0$. Hence we may again define a $G$-invariant positive-definite Hermitian form on $K[x_1, \ldots, x_n, \theta_1, \ldots, \theta_n]$ by
  \[ (f, \omega) \coloneqq \text{constant coefficient of }\partial_\omega f. \]
This form is linear in $f$ and conjugate-linear in $\omega$.

Under this form, the super coinvariant ideals and super harmonics from \Cref{ssec:intro_super_harmonics} are orthogonal complements, $\cJ_G^\perp = \cSH_G$ and $\cSH_G^\perp = \cJ_G$. Furthermore, the adjoint of $\partial_{\theta_i}$ with respect to this form is left multiplication by $\theta_i$, i.e. $\partial_{\theta_i}^\dagger f = \theta_i f$. More generally, $\partial_\omega^\dagger = \omega f$. See \cite[Lemmas~2.7, 2.11, 5.4]{MR4258767} for more information.

\subsection{Generalized exterior derivatives}\label{ssec:generalized_exterior_derivatives}

Recall the exterior derivative \eqref{eq:exterior_derivative}, which is defined by $\dif = \sum_{i=1}^n \partial_{x_i} \theta_i$. Here and elsewhere, $\theta_i$ refers to left multiplication by $\theta_i$. Equation \eqref{eq:exterior_derivative} is consistent with our  usage of $\dif x_i$ in \Cref{ssec:differential_forms}. We now describe the operators $\dif_i$ generalizing \eqref{eq:exterior_derivative} and underlying our construction of $\cSH_G^{\det}$ in \cite{MR4258767}.

We have $\sigma \cdot \dif f = \dif(\sigma \cdot f)$ for all $f \in K[x_1, \ldots, x_n, \theta_1, \ldots, \theta_n]$. The $G$-equivariance of the exterior derivative conceptually arises from the $G$-invariance of the form $\overline{x}_1\theta_1 + \cdots + \overline{x}_n\theta_n$, or equivalently the $G$-invariance of the norm-squared function $\overline{x}_1 x_1 + \cdots + \overline{x}_n x_n$, which is $G$-invariant since $G \subset U(n, K)$ consists of unitary matrices. Here $\overline{x}_i(v) \coloneqq \overline{x_i(v)}$.

More generally, every $G$-invariant element of $K[\overline{x}_1, \ldots, \overline{x}_n, \theta_1, \ldots, \theta_n]$ gives rise to a corresponding $G$-equivariant operator in
  \[ K[\partial_{x_1}, \ldots, \partial_{x_n}, \theta_1, \ldots, \theta_n] \]
obtained by replacing $\overline{x}_i$ with $\partial_{x_i}$. We now summarize the structure of these operators and define the operators $\dif_i$ used in \Cref{ssec:intro_semi-invariants}.

When $G$ is a pseudo-reflection group, Shephard--Todd \cite{MR0059914} and later Chevalley \cite{MR72877} showed that $K[x_1, \ldots, x_n]^G = K[f_1, \ldots, f_n]$ for algebraically independent homogeneous $G$-invariants $f_1, \ldots, f_n$ called \textit{basic invariants} of $G$. The basic invariants are not unique, though the multiset $\{d_1, \ldots, d_n\}$ of their degrees is uniquely determined and is called the multiset of \textit{degrees} of $G$. The \textit{exponents} of $G$ are the multiset $\{d_i - 1 , \ldots, d_n - 1\}$, which is the multiset of the $x$-degrees of $\dif f_i$.

Solomon \cite{MR0154929} described the $G$-invariants $K[x_1, \ldots, x_n, \theta_1, \ldots, \theta_n]^G$. He showed that they have a $K$-basis given by
\begin{equation}\label{eq:Solomon_invariants}
  \{f_1^{\alpha_1} \cdots f_n^{\alpha_n} \dif f_{i_1} \cdots \dif f_{i_k} : \alpha_i \geq 0, 1 \leq i_1 < \cdots < i_k \leq n, k \in [n]\}.
\end{equation}
Consequently,
\begin{equation}\label{eq:Solomon_generators}
  \cJ_G = \langle K[x_1, \ldots, x_n, \theta_1, \ldots, \theta_n]_+^G\rangle = \langle f_1, \ldots, f_n, \dif f_1, \ldots, \dif f_n\rangle.
\end{equation}

Orlik--Solomon \cite{MR575083} generalized Solomon's exterior algebra construction of the invariants to certain Galois conjugates. Translated to the present language of differential operators, we have the following.

\begin{theorem}[{\cite{MR575083}}; see {\cite[\S3.1]{MR4258767}}]\label{thm:Orlik-Solomon}
  There are bi-homogeneous, $G$-equivariant operators $\dif_1, \ldots, \dif_n$ which raise $\theta$-degree by $1$ such that a $K$-basis for the $G$-equivariant operators in $K[\partial_{x_1}, \ldots, \partial_{x_n}, \theta_1, \ldots, \theta_n]^G$ is given by
    \[ \{\partial_{f_1}^{\alpha_1} \cdots \partial_{f_n}^{\alpha_n} \dif_{i_1} \cdots \dif_{i_k} : \alpha_i \geq 0, 1 \leq i_1 < \cdots < i_k \leq n, k \in [n]\}. \]
\end{theorem}

\begin{definition}
  The $\dif_i$ appearing in \Cref{thm:Orlik-Solomon} will be called \textit{generalized exterior derivative} on $K[x_1, \ldots, x_n, \theta_1, \ldots, \theta_n]$.
\end{definition}

See \Cref{tab:explicit} for an explicit description of the generalized exterior derivatives for $G(m, p, n)$. For example, when $G = \fS_n$, one may use
  \[ \dif_i = \sum_{j=1}^n \partial_{x_j}^i \theta_j \]
where $0 \leq i \leq n-1$; see \Cref{sec:background:explicit} for further details. When $K \subset \mathbb{R}$, so $G$ is a (real) reflection group, then by \eqref{eq:Solomon_invariants} we may take $\dif_i = \sum_{j=1}^n \partial_{f_{ij}} \theta_j$ where $f_{ij} = \partial_{x_j} f_i$.

Since the generalized exterior derivatives raise $\theta$-degree by $1$, they satisfy $\dif_i \dif_j = -\dif_j \dif_i$, so in particular $\dif_i^2 = 0$. The $\dif_i$ are not unique, though the multiset of degrees $\{e_1^*, \ldots, e_n^*\}$ by which they lower $x$-degree is uniquely determined. This multiset by definition consists of the \textit{co-exponents} of $G$.

\subsection{Removing invariant vectors}\label{ssec:intro_removing_invariant_vectors}

Our description of the semi-invariant differential forms in \Cref{thm:semi-invariants} from \Cref{ssec:intro_semi-invariants} involves $r$ operators $\dif_1, \ldots, \dif_r$ where $r = \dim(V/V^G) \leq n$, rather than $n$ operators. This arises from a slight mismatch between \Cref{thm:semi-invariants} and the result \cite[Thm.~5.7]{MR4258767} underlying it. That result involves semi-invariant harmonic forms in $\rS(V^*) \otimes \wedge M^*$, where $M$ is a finite-dimensional $G$-module with $M^G = 0$. The generalization of the basis \eqref{eq:SH_det_basis} from \cite{MR4258767} involves $\dim M$ differential operators in general. When $G = \fS_n$ acts on $\mathbb{Q}^n$ by permutation matrices, $V^G = \Span_{\bQ}\{(1, \ldots, 1)\} \neq 0$, so the result does not directly apply. However, we may use $M = V/V^G$ in \cite[Thm.~5.7]{MR4258767}, which in the case of $G = \fS_n$ is the \textit{standard representation}.

The relationship between the coinvariants and harmonics of $\rS(V^*) \otimes \wedge (V/V^G)^*$ and $\rS(V^*) \otimes \wedge V^*$ is straightforward. Write $\cJ_G'$, $\cSH_G'$, and $\cSR_G'$ for the super coinvariant ideal, super harmonics, and super coinvariant algebra of $\cS(V^*) \otimes \wedge (V/V^G)^*$.

\begin{lemma}\label{lem:removing_invariant_vectors}
  The super harmonics $\cSH_G$ and $\cSH_G'$ are naturally isomorphic.
  
  \begin{proof}(Sketch.)
    We have a natural map $\Phi \colon \rS(V^*) \otimes \wedge V^* \to \rS(V^*) \otimes \wedge (V/V^G)^*$. By replacing $\theta_1, \ldots, \theta_n$ with a basis $\psi_1, \ldots, \psi_r, \psi_{r+1}, \ldots, \psi_n$ for which $\psi_{r+1}, \ldots, \psi_n$ is a basis for $\wedge^1 (V^G)^*$, we may think of $\rS(V^*) \otimes \wedge (V/V^G)^*$ as being obtained from $\rS(V^*) \otimes \wedge V^*$ by setting $\psi_1, \ldots, \psi_r = 0$. Moreover, we may suppose $\Span_K\{\psi_1, \ldots, \psi_r\}$ is $G$-stable. We then find that $\Phi(\cJ_G) = \cJ_G'$, $\Phi$ descends to an isomorphism $\cSR_G \too{\sim} \cSR_G'$, and $\Phi$ restricts to an isomorphism $\cSH_G \too{\sim} \cSH_G'$.
  \end{proof}
\end{lemma}

In particular, the $\det$-isotypic component $\cSH_G^{\det}$ has dimension $2^r$ rather than $2^n$. Similarly, the ``volume form'' in $\cSH_G$ is an $r$-form rather than an $n$-form, namely $\psi_1 \cdots \psi_r$ in the notation of the proof, and $\cSH_G^r$ is the top non-zero component of $\cSH_G$.

We may hence choose generalized exterior derivatives $\dif_{r+1}, \ldots, \dif_n$ which each fix the $x$-degree and act as $0$ on $\cSH_G$. The corresponding co-exponents $e_{r+1}^*, \ldots, e_n^*$ are all zero and will be ignored. The remaining generalized exterior derivatives $\dif_1, \ldots, \dif_r$ strictly decrease $x$-degree and their co-exponents are positive. The description of $\cSH_G^{\det}$ in \eqref{eq:SH_det_basis} now follows from \cite[Thm.~5.7]{MR4258767}.

\subsection{Vandermondians and Jacobians}\label{ssec:Vandermondians}

We now briefly give an explicit construction of the key element $\Delta_G$ from \Cref{ssec:intro_harmonics}, which we call the \textit{Vandermondian} of $G$, along with a related element $\Delta_G^*$, which we call the \textit{co-Vandermondian} of $G$. See \cite[\S3]{MR4258767} for a more complete summary and references to the literature. We continue the notation from \Cref{ssec:generalized_exterior_derivatives}, so $G$ is a pseudo-reflection group and $f_1, \ldots, f_n$ are basic invariants of $G$.

By Steinberg's \Cref{thm:Steinberg}, the top-degree component of $\cH_G \subset K[x_1, \ldots, x_n]$ is spanned by an element $\Delta_G \in \cH_G$, uniquely defined up to non-zero scalar multiples, which transforms according to $g \cdot \Delta_G = \det(g) \Delta_G$. Similarly, there is an element $\Delta_G^* \in \cH_G$, uniquely defined up to non-zero scalar multiples, which transforms according to $g \cdot \Delta_G^* = \overline{\det(g)} \Delta_G$. By Steinberg's result, $\deg \Delta_G \geq \deg \Delta_G^*$. In fact, $\Delta_G^* \mid \Delta_G$ and equality holds if and only if $G$ is generated by reflections (that is, order $2$ pseudo-reflections). These facts may be read off from a formula of Gutkin \cite{MR0314956}, which expresses $\Delta_G$ and $\Delta_G^*$ explicitly in terms of the reflecting hyperplanes of $G$ as follows.

Let $\cA(G)$ be the set of reflecting hyperplanes of $G$, i.e.~the fixed spaces of pseudo-reflections of $G$. For each $H \in \cA(G)$, fix some $\alpha_H \in V^*$ with $\ker \alpha_H = H$. Let $G_H$ denote the subgroup of $G$ fixing $H$ pointwise. The Vandermondian is defined uniquely up to a non-zero scalar by
  \[ \Delta_G = \prod_{H \in \cA(G)} \alpha_H^{|G_H| - 1}. \]
The co-Vandermondian is defined uniquely up to a non-zero scalar by
  \[ \Delta_G^* = \prod_{H \in \cA(G)} \alpha_H. \]

Furthermore, $\deg \Delta_G = \sum_{i=1}^n (\deg(f_i) - 1)$ is the sum of the exponents of $G$, which is the $x$-degree of $\dif f_1 \cdots \dif f_n$. Similarly, the sum of the co-exponents is the number of reflecting hyperplanes, $e_1^* + \cdots + e_n^* = |\cA(G)|$, which is the amount the $x$-degree is lowered by $\dif_1 \cdots \dif_r$.

Given a set of homogeneous $G$-invariants $f_1, \ldots, f_n$, one may verify that they are indeed basic invariants using Saito's criterion \cite[Thm.~4.19]{MR1217488} or an appropriate generalization \cite[Thm.~3.1]{MR575083}, which says that it suffices to check that the Jacobian determinant of $f_1, \ldots, f_n$ agrees with the Vandermondian of $G$. That is, we require $J_G \coloneqq \det(\partial_{x_j} f_i)_{1 \leq i, j \leq n} = \Delta_G$, up to a non-zero constant. Note that in this case $\dif f_1 \cdots \dif f_n$ is a non-zero multiple of $\Delta_G \,\theta_1 \cdots \theta_n$.

Likewise, for the generalized exterior derivatives $\dif_i = \sum_{j=1}^n \partial_{f_{ij}} \theta_j$, the corresponding Jacobian determinant $J_G^* \coloneqq \det(f_{ij})_{1 \leq i, j \leq n}$ transforms according to $\sigma \cdot J_G^* = \overline{\det(\sigma)} J_G^*$. Given a proposed set of bi-homogeneous $G$-equivariant operators, one may verify that they are indeed generalized exterior derivatives by checking that $J_G^*$ agrees with the co-Vandermondian of $G$, $\det(f_{ij})_{1 \leq i, j \leq n} = \Delta_G^*$, up to a non-zero constant. Note that $\dif_1 \cdots \dif_n = \partial_{\Delta_G^*} \,\theta_1 \cdots \theta_n$.

\subsection{Explicit formulas}\label{sec:background:explicit}

We now give explicit descriptions for the basic invariants, generalized exterior derivatives, Vandermondians, co-Vandermondians, and co-exponents of the pseudo-reflection groups $G = G(m, p, n)$ described in \Cref{ssec:intro:psgs}. See \Cref{tab:explicit} at the end of the paper for a quick summary. The special cases of real reflection groups $G(1, 1, n) = \fS_n$, $G(2, 1, n) = \fB_n$, $G(2, 2, n) = \fD_n$, and $G(m, m, 2) = \fDih_m$ are also written out in \Cref{tab:explicit}.

When $G = G(1, 1, n) = \fS_n$ is the symmetric group, we may use power-sums for the basic invariants, namely $f_i = \sum_{j=1}^n x_j^i$ for $1 \leq i \leq n$. The reflections are the $\binom{n}{2}$ transpositions with reflecting hyperplanes $x_j - x_i = 0$, so $\Delta_{\fS_n} = \prod_{1 \leq i < j \leq n} (x_j - x_i)$ is the classical Vandermonde determinant of degree $\binom{n}{2}$. Since $G$ is a real reflection group, the co-Vandermondian is equal to the Vandermondian, and we may use the exterior derivatives of the $f_i$ to construct the generalized exterior derivatives $\dif_i$. It is most convenient to use $\dif_i = \sum_{j=1}^n \partial_{x_j^i} \theta_j$ for $0 \leq i \leq n-1$. In particular, $\dif_1 = \dif$ is the exterior derivative and $\dif_0 = \theta_1 + \cdots + \theta_n$ reflects the fact that $x_1 + \cdots + x_n$ is $\fS_n$-invariant. Note that $\dif_0$ acts as $0$ on $\cSH_{\fS_n}$, so only $\dif_1, \ldots, \dif_{n-1}$ are of interest, and $r = n-1$. The co-exponents are $e_i^* = i$ and their sum is $\binom{n}{2}$.

When $G = G(m, 1, n)$ is the group of pseudo-permutation matrices whose non-zero entries belong to the group $\mu_m$ of $m$th complex roots of unity, we may use basic invariants $f_i = \sum_{j=1}^n x_j^{mi}$ for $1 \leq i \leq n$. The pseudo-reflections come in two types. First, the $m\binom{n}{2}$ generalized transpositions indexed by $1 \leq i < j \leq n$ and $a \in \mu_m$ where $\sigma(f_i) = af_j, \sigma(f_j) = a^{-1}f_i, \sigma(f_k) = f_k$ for $k \not\in \{i, j\}$, which have reflecting hyperplanes $ax_i - x_j$ with $|G_H| = 2$. Second, the $n(m-1)$ ``rotations'' indexed by $1 \leq i \leq n$ and $a \in \mu_m - \{1\}$ where $\sigma(f_i) = af_i$, $\sigma(f_k) = f_k$ for $k \neq i$, which have reflecting hyperplanes $x_i = 0$ with $|G_H| = m$. Now $\Delta_{G(m, 1, n)} = (x_1 \cdots x_n)^{m-1}\prod_{1 \leq i < j \leq n} (x_j^m - x_i^m)$, which has degree $m\binom{n}{2} + n(m-1)$. When $m > 1$, the co-Vandermondian is $\Delta_{G(m, 1, n)}^* = (x_1 \cdots x_n) \prod_{1 \leq i < j \leq n} (x_j^m - x_i^m)$, which has degree $m \binom{n}{2} + n$. In this case, the generalized exterior derivatives are $\dif_i = \sum_{j=1}^n \partial_{x_j^{(i-1)m+1}} \theta_j$ for $1 \leq i \leq n$, and the co-exponents are $e_i^* = (i-1)m+1$.

When $G = G(m, p, n)$ is the subgroup of $G(m, 1, n)$ where the product of the non-zero entries raised to the $(m/p)$th power is $1$, we may use basic invariants $f_i = \sum_{j=1}^n x_j^{mi}$ for $1 \leq i \leq n-1$ and $f_n = (x_1 \cdots x_n)^{m/p}$. The pseudo-reflections come in the same types as for $G(m, 1, n)$, except that there are $n(m/p - 1)$ ``rotations'' which additionally require $a \in \mu_{m/p} - \{1\}$, and $|G_H| = m/p$. Now $\Delta_{G(m, p, n)} = (x_1 \cdots x_n)^{m/p - 1} \prod_{1 \leq i < j \leq n} (x_j^m - x_i^m)$, which has degree $m\binom{n}{2} + n(m/p - 1)$. When $p \neq m$, the co-Vandermondian is $x_1 \cdots x_n \prod_{1 \leq i < j \leq n} (x_j^m - x_i^m)$, which has degree $m \binom{n}{2} + n$. In this case, the generalized exterior derivatives are $\dif_i = \sum_{j=1}^n \partial_{x_j^{(i-1)m+1}} \theta_j$, and the co-exponents are $e_i^* = (i-1)m+1$. When $p=m$, the only pseudo-reflections are the generalized transpositions and the co-Vandermondian is $\prod_{1 \leq i < j \leq m} (x_j^m - x_i^m)$, which has degree $m \binom{n}{2}$. In this case, the generalized exterior derivatives are $\dif_i = \sum_{j=1}^n \partial_{x_j^{(i-1)m+1}} \theta_j$ for $1 \leq i \leq n-1$ and $\dif_n = \sum_{j=1}^n \partial_{(x_1 \cdots \widehat{x}_j \cdots x_n)^{m-1}} \theta_j$, with co-exponents $e_i^* = (i-1)m+1$ for $1 \leq i \leq n-1$ and $e_n^* = (n-1)(m-1)$.

\section{The structure of $\cSH_G^r$}\label{sec:top}

We begin by considering the top $\theta$-degree component of the key differential operator equation \eqref{eq:thm:A}. Let $G$ be a pseudo-reflection group with basic invariants $f_1, \ldots, f_n$. By \Cref{lem:removing_invariant_vectors}, we may streamline our exposition by supposing without loss of generality throughout this subsection that $V^G = 0$. Hence $r=n$ and the basic invariants $f_1, \ldots, f_n$ all have degree at least $2$.

Recall that $\Delta_G, \Delta_G^* \in \cH_G \subset K[x_1, \ldots, x_n]$ are the unique elements up to non-zero scalar multiples in $\cH_G^{\det}$ and $\cH_G^{\overline{\det}}$, respectively. Since $\cH_G$ is isomorphic to the regular representation of $G$, there is also a non-zero element
  \[ \Gamma_G \in \cH_G^{\det^2} \]
unique up to non-zero scalar multiples. We may describe $\Gamma_G$ in terms of $\Delta_G$ and $\Delta_G^*$ as follows.

\begin{lemma}\label{lem:Gamma_G}
  We have $\Gamma_G = \partial_{\Delta_G^*} \Delta_G$ up to non-zero scalar multiples.
  
  \begin{proof}
    Since $\Delta_G \in \cH_G$, we have $\partial_{\Delta_G^*} \Delta_G \in \cH_G$. We have
    \begin{align*}
      \sigma \cdot \partial_{\Delta_G^*} \Delta_G
        &= \partial_{\sigma(\Delta_G^*)} \sigma(\Delta_G) \\
        &= \partial_{\overline{\det(\sigma)} \Delta_G^*} \det(\sigma) \Delta_G \\
        &= \det(\sigma)^2 \partial_{\Delta_G^*} \Delta_G,
    \end{align*}
    so $\partial_{\Delta_G^*} \Delta_G \in \cH_G^{\det^2}$. By Gutkin's formula in \Cref{ssec:Vandermondians}, $\Delta_G^* \mid \Delta_G$, so
      \[ 0 \neq \partial_{\Delta_G} \Delta_G = \partial_{\Delta_G/\Delta_G^*} \partial_{\Delta_G^*} \Delta_G. \]
    Hence $\partial_{\Delta_G^*} \Delta_G \neq 0$.
  \end{proof}
\end{lemma}

Consider the following analogues of the classical coinvariant ideal $\cI_G$ and harmonics $\cH_G$.

\begin{definition}
  Let
    \[ \cI_G' \coloneqq \langle \partial_{x_j} f_i : i, j \in [n]\rangle \]
  and
    \[ \cH_G' \coloneqq \{f \in K[x_1, \ldots, x_n] : \partial_g f = 0 \text{ for all }g \in \cI_G'\}. \]
\end{definition}

\noindent The ideal $\cI_G'$ is the $(n-1)$st Fitting ideal of the Jacobian of the basic invariants $f_1, \ldots, f_n$. Since $\deg f_i \geq 2$, $\cI_G'$ is proper. It is independent of the basis $x_1, \ldots, x_n$.

\begin{lemma}\label{lem:H_G_prime_props}
  We have $\cI_G' \supset \cI_G$, $\cH_G' \subset \cH_G$, and $K[x_1, \ldots, x_n] = \cH_G' \oplus \cI_G'$.
  
  \begin{proof}
    Recall Euler's formula, $\sum_{j=1}^n x_j \partial_{x_j} f = \deg(f) f$ for homogeneous $f$. Thus $f_i \in \cI_G'$ for all $i$, so $\cI_G' \supset \cI_G$. We now see directly that $\cH_G' \subset \cH_G$. Finally, $(\cI_G')^\perp = \cH_G'$.
  \end{proof}
\end{lemma}

\begin{lemma}\label{lem:top_harmonics}
  Suppose $V^G = 0$. Then
    \[ \cSH_G^n = \cH_G' \theta_1 \cdots \theta_n. \]
  
  \begin{proof}
    The $n$-forms $\omega \in \cSH_G^n$ are of the form $g\,\theta_1 \cdots \theta_n$ for some $g \in K[x_1, \ldots, x_n]$. We have $\omega \in \cSH_G^n$ if and only if $\partial_{f_i} \omega = 0$ and $\partial_{\dif f_i} \omega = 0$ for all $1 \leq i \leq n$. The first condition occurs if and only if $g \in \cH_G$. Noting that $\dif f_i = \sum_{j=1}^n \partial_{x_j} f_i \theta_j$, the second condition implies that
    \begin{align*}
      \partial_{\dif f_i} g \,\theta_1 \cdots \theta_n
         &= \sum_{j=1}^n \pm \partial_{\partial_{x_j} f_i} g \,\theta_1 \cdots \widehat{\theta}_j \cdots \theta_n.
    \end{align*}
    Hence $\partial_{\dif f_i} \omega = 0$ if and only if $\partial_{\partial_{x_j} f_i} g = 0$ for all $j$, so the second condition is equivalent to $g \in \cH_G'$. The result follows since $\cH_G' \subset \cH_G$.
  \end{proof}
\end{lemma}

\begin{definition}
  For $F \in K[x_1, \ldots, x_n]$, write
    \[ \Ann F \coloneqq \{g \in K[x_1, \ldots, x_n] : \partial_g F = 0\}. \]
\end{definition}

Steinberg proved \Cref{thm:Steinberg} by first showing \cite[Thm.~1.3(b)]{MR167535}
\begin{equation}\label{eq:Steinberg_ann}
  \Ann \Delta_G = \cI_G.
\end{equation}
We now consider $\Ann \Gamma_G$.

\begin{lemma}\label{lem:Gamma_G_in_Hprime}
  Suppose $V^G = 0$. Then $\Gamma_G \in \cH_G'$. Equivalently,
    \[ \Ann \Gamma_G \supset \cI_G'. \]
  
  \begin{proof}
    As noted in \Cref{ssec:Vandermondians}, $\dif_1 \cdots \dif_n = \partial_{\Delta_G^*} \theta_1 \cdots \theta_n$, so $\Gamma_G \theta_1 \cdots \theta_n = \partial_{\Delta_G^*} \Delta_G \theta_1 \cdots \theta_n = \dif_1 \cdots \dif_n \Delta_G$. By \eqref{eq:SH_det_basis} and \Cref{lem:top_harmonics}, $\dif_1 \cdots \dif_n \Delta_G \in \cSH_G^n = \cH_G' \theta_1 \cdots \theta_n$, so $\Gamma_G \in \cH_G'$. Hence $\partial_g \Gamma_G = 0$ for all $g \in \cI_G'$, so $\cI_G' \subset \Ann \Gamma_G$.
  \end{proof}
\end{lemma}

We will now show that $\Ann \Gamma_G = \cI_G'$ is equivalent to the $k=n$ case of \eqref{eq:thm:A} when $V^G = 0$. We then restate and prove \Cref{thm:A2}.

\begin{proposition}\label{prop:H_prime_annihilator}
  Suppose that $G$ is a pseudo-reflection group and $V^G = 0$. Then
  \begin{equation}\label{eq:H_prime_annihilator.0}
    \cSH_G^n = K[\partial_{x_1}, \ldots, \partial_{x_n}] (\cSH_G^n)^{\det}
  \end{equation}
  if and only if
  \begin{equation}\label{eq:H_prime_annihilator}
    \Ann \Gamma_G = \cI_G'.
  \end{equation}
  
  \begin{proof}
    By \eqref{eq:SH_det_basis} and \Cref{lem:Gamma_G},
    \begin{align*}
     (\cSH_G^n)^{\det}
       &=K\,\dif_1 \cdots \dif_n \Delta_G \\
       &= K\,\partial_{\Delta_G^*} \Delta_G \theta_1 \cdots \theta_n \\
       &= K\,\Gamma_G \theta_1 \cdots \theta_n.
    \end{align*}
    Hence
       \[ K[\partial_{x_1}, \ldots, \partial_{x_n}] (\cSH_G^n)^{\det}
           = K[\partial_{x_1}, \ldots, \partial_{x_n}] \Gamma_G \theta_1 \cdots \theta_n. \]
    By \Cref{lem:top_harmonics}, equation \eqref{eq:H_prime_annihilator.0} is hence equivalent to
    \begin{equation}\label{eq:H_prime_annihilator.2}
      \cH_G' = \{\partial_g \Gamma_G : g \in K[x_1, \ldots, x_n]\}.
    \end{equation}
    
    Consider the conjugate-linear map $\psi \colon K[x_1, \ldots, x_n] \to K[x_1, \ldots, x_n]$ given by $\psi(g) \coloneqq \partial_g \Gamma_G$, so $\ker \psi = \Ann \Gamma_G$. Equation \eqref{eq:H_prime_annihilator.2} is equivalent to $\im\psi = \cH_G'$. By \Cref{lem:Gamma_G_in_Hprime}, we have $\psi \colon K[x_1, \ldots, x_n] \to \cH_G'$. Again by \Cref{lem:Gamma_G_in_Hprime}, $\ker\psi \supset \cI_G'$. Thus by \Cref{lem:H_G_prime_props}, $\im\psi = \im\psi|_{\cH_G'}$ and $\ker\psi = \ker\psi|_{\cH_G'} \oplus \cI_G'$. Since $\cH_G'$ is finite-dimensional, $\psi|_{\cH_G'}$ is surjective if and only if it is injective, which occurs if and only if $\ker\psi = \cI_G'$. This is a restatement of \eqref{eq:H_prime_annihilator}.
  \end{proof}
\end{proposition}

\begin{reptheorem}{thm:A2}
  Let $G = G(m, 1, n)$ or let $G$ be real. Then the $\theta$-degree $r$ component of \eqref{eq:thm:A} holds.
  
  \begin{proof}
    First suppose $G$ is real. Then $\det^2 = 1$, so $\Gamma_G = 1$ and $\Ann \Gamma_G = \langle x_1, \ldots, x_n\rangle$. Furthermore, $x_1^2 + \cdots + x_n^2$ is $G$-invariant, so $\cI_G' = \langle x_1, \ldots, x_n\rangle$ as well. The result follows from \Cref{prop:H_prime_annihilator}.
    
    If $G = G(m, 1, n)$ for $m > 1$, then $r=n$. \Cref{tab:explicit} gives $\cI_G' = \langle x_1^{m-1}, \ldots, x_n^{m-1}\rangle$. Furthermore, we claim $\Gamma_G = (x_1 \cdots x_n)^{m-2}$. First, this is annihilated by $\partial_{f_i}$ by degree considerations, so it is harmonic. Also, it transforms under $G$ by $\det^2$. Indeed, if $\sigma(f_i) = \lambda_i f_{\tau(i)}$, $\det(\sigma) = \lambda_1 \cdots \lambda_n \det(\tau)$, so $\det(\sigma)^2 = (\lambda_1 \cdots \lambda_n)^2$ and
      \[ \sigma (x_1 \cdots x_n)^{m-2} = (\overline{\lambda_1} \cdots \overline{\lambda_n} x_1 \cdots x_n)^{m-2} = (\lambda_1 \cdots \lambda_n)^2 (x_1 \cdots x_n)^{m-2}. \]
  Finally, the annihilator of $(x_1 \cdots x_n)^{m-2}$ is precisely $\langle x_1^{m-1}, \ldots, x_n^{m-1}\rangle$.
  \end{proof}
\end{reptheorem}

The assertion is false for $G(m, p, n)$ not covered by the above theorem. Here we consider the cyclic groups $G(m, p, 1) = G(m/p, 1, 1)$ as part of the family $G(m, 1, n)$. The simplest example with strict containment is $G(4, 2, 2)$, which is generated by reflections but is not real.

\begin{lemma}\label{lem:no_dice}
  If $G = G(m, p, n)$ is not of the form $G(1, 1, n)$, $G(2, 1, n)$, $G(2, 2, n)$, $G(m, p, 1) = G(m/p, 1, 1)$, $G(m, m, 2)$, or $G(m, 1, n)$, then $r=n$ and $\Ann \Gamma_G \neq \cI_G'$. Hence \eqref{eq:H_prime_annihilator} does not hold, so the $r$-form component of \eqref{eq:thm:A} does not hold.
  
  Moreover, in these cases, the top $x$-degree component of $\cSH_G^r$ is strictly higher than the top $x$-degree component of $(\cSH_G^r)^{\det}$, so the $r$-form component of \eqref{eq:lem:B} does not hold.
  
  \begin{proof}
    By \Cref{tab:explicit}, we may use basic invariants $f_i = \sum_{j=1}^n x_j^{mi}$ for $1 \leq i \leq n-1$ and $f_n = (x_1 \cdots x_n)^{m/p}$. Hence
      \[ \cI_G' = \langle x_1^{m-1}, \ldots, x_n^{m-1}, x_1^{-1} (x_1 \cdots x_n)^{m/p}, \ldots, x_n^{-1} (x_1 \cdots x_n)^{m/p}\rangle. \]
    
    If $p=m$, $\Delta_G = \Delta_G^*$ (even though $G(m, m, n)$ is typically not real), so $\Gamma_G = 1$ and $\Ann \Gamma_G = \langle x_1, \ldots, x_n\rangle$. The top $x$-degree component of $(\cSH_G^n)^{\det}$ is hence degree $0$. By assumption, $m \geq 3$, so the generators $x_i^{m-1}$ of $\cI_G'$ are at least quadratic. Also by assumption, $n \geq 3$, so the generators $x_i^{-1} (x_1 \cdots x_n)$ are also at least quadratic, so $\Ann \Gamma_G \supsetneq \cI_G'$. Moreover, $\cI_G'$ does not contain the linear polynomials, so $\cH_G'$ contains all linear polynomials, and the top $x$-degree component of $\cSH_G^n$ is at least degree $1$.
    
    If $p \neq m$, one finds $\Gamma_G = (x_1 \cdots x_n)^{m/p - 2}$ as in the proof of \Cref{thm:A2}, so $\Ann \Gamma_G = \langle x_1^{m/p-1}, \ldots, x_n^{m/p-1}\rangle$. The generators of $\cI_G'$ have strictly larger degree since $n \geq 2$, so the containment is strict. Indeed, it is easy to see that $x_1 (x_1 \cdots x_n)^{m/p - 2} \not\in \cI_G'$, so the top $x$-degree component of $\cSH_G^n$ is at least $1$ higher than that of $(\cSH_G^n)^{\det}$.
  \end{proof}
\end{lemma}

\begin{example}
  When $G = G(4, 2, 2)$, we have $\cI_G' = \langle x_1^3, x_2^3, x_1^2 x_2, x_1 x_2^2 \rangle$ while $\Gamma_G = 1$ has annihilator $\langle x_1, x_2\rangle$. Hence $\cH_G'$ consists of all $g$ of degree at most $2$ and $\cSH_G^2 = \{g\,\theta_1\theta_2 : \deg g \leq 2\}$. On the other hand, $(\cSH_G^2)^{\det} = K\,\theta_1 \theta_2$. Hence \eqref{eq:thm:A} and \eqref{eq:lem:B} are false in this case.
\end{example}

We may also use the explicit description of $\cH_G'$ to determine the highest degree of multiples of the volume form in $\cSH_G$. This will be used below in \Cref{sec:total}. 

\begin{lemma}\label{lem:top_xdeg}
  Let $G = G(m, p, n)$. Then the top $x$-degree component of $\cSH_G^r$ has degree
  \begin{align*}
    \max \{i : \cSH_G^{i, r} \neq 0\} = 
    \begin{cases}
      0 & \text{if }n=1, m/p \leq 2 \\
      m/p-2 & \text{if }n=1, m/p \geq 3 \\
      0 & \text{if }n \geq 2, m \leq 2 \\
      m/p-2 + (n-1)(m-2) & \text{if }n \geq 2, m \geq 3.
    \end{cases}
  \end{align*}
  
  \begin{proof}
    If $n = 1 = m/p$, then $\cSH_G^0 = K$ has degree $0$. If $n=1$ and $m/p \geq 2$, then $f_1 = x_1^{m/p}$ with $\dif f_1 = m/p \cdot x_1^{m/p-1} \theta_1$. Hence $x_1^s \in \cSH_G^1$ if and only if $\partial_{x_1}^{m/p-1} x^s = 0$, or if and only if $s \leq m/p-2$. 
    
    Now take $n \geq 2$. As above, we have
      \[ \cI_G' = \langle x_1^{m-1}, \ldots, x_n^{m-1}, x_1^{-1} (x_1 \cdots x_n)^{m/p}, \ldots, x_n^{-1} (x_1 \cdots x_n)^{m/p} \rangle \]
    for $(m, p) \neq (1, 1)$. If $m=p=1$, then $G$ is real and $\cSH_G^r = K$ has top degree $0$ by \Cref{thm:A2}. If $p=1$ and $m \geq 2$, we have $\cI_G' = \langle x_1^{m-1}, \ldots, x_n^{m-1}\rangle$, so the top-degree monomial not in $\cI_G'$ is $(x_1 \cdots x_n)^{m-2}$ which has degree $n(m-2)$. If $p=m=2$, then $\cI_G' = \langle x_1, \ldots, x_n\rangle$ and again $\cSH_G^r = K$.
    
    Now take $n \geq 2$, $m \geq 3$, and $p \geq 2$. Here $m/p \leq m/2 < m-1$. Suppose $x^\alpha \not\in \cI_G'$ for $|\alpha|$ maximal. By symmetry, we may suppose $\alpha_1 \leq \cdots \leq \alpha_n$. If $\alpha_1 \geq m/p$, then $(x_1 \cdots x_n)^{m/p} \mid x^\alpha$ and $x^\alpha = 0$, so $\alpha_1 \leq m/p-1$. If $\alpha_1 = m/p-1$, then $\alpha_2 = m/p-1$, and we find $x^\alpha = x_1^{m/p-1} x_2^{m/p-1} (x_3 \cdots x_n)^{m-2}$. If $\alpha_1 = m/p-2$, then we find $x^\alpha = x_1^{m/p-2} (x_2 \cdots x_n)^{m-2}$. The degree of the latter minus the degree of the former is $m - 2 - m/p \geq 0$, so the latter is the top-degree element.
  \end{proof}
\end{lemma}

\begin{corollary}\label{cor:top_xdeg}
  Let $G = G(m, p, n)$. The top total degree component of $\cSH_G^r$ has degree strictly below $\deg \Delta_G$, except when $n=1$ or $(m, p, n) \in \{(1, 1, 2), (2, 2, 2)\}$ when equality holds.
  
  \begin{proof}
    In this case,
      \[ \deg \Delta_G = m \binom{n}{2} + n\left(\frac{m}{p} - 1\right). \]
    First take $n=1$. If $m=p=1$, then $r=0$ and both degrees are $0$. If $m/p \geq 2$, then $r=1$, $\deg \Delta_G = m/p-1$, and the top \textit{total} degree of $\cSH_G^r$ is $m/p-2+1 = m/p-1$ by \Cref{lem:top_xdeg}, so equality holds.
    
    Now consider $n \geq 2$. If $m=1$, then $\binom{n}{2} \geq r = n-1$, and the inequality is strict for $n \geq 3$. If $m=2$, then $m \binom{n}{2} + n(m/p-1) \geq n$, and the inequality is strict except when $m=p=n=2$. If $m \geq 3$, we have
    \begin{align*}
      \left(m \binom{n}{2} + n\left(\frac{m}{p} - 1\right)\right) &- \left(\frac{m}{p} - 2 + (n-1)(m-2) + n\right)\\
        &= \frac{m (n-1) (p(n-2)+2)}{2p} > 0.
    \end{align*}
  \end{proof}
\end{corollary}

\section{Exactness of exterior differentiation on $\cSH_G$}\label{sec:complex}

In this section we introduce a generalization of the complex \eqref{eq:exact_d} for a pseudo-reflection group $G$. We then summarize an algebraic analogue of well-known results in Hodge theory. Finally we prove \Cref{thm:exact} and deduce \Cref{thm:A1}. Our argument follows an approach to proving exactness of the Koszul complex; see for example \cite[Prop.~7.2.11]{MR1606831}.

\subsection{Super harmonic cochain complexes}

We begin by showing that two ``dual'' types of operators preserve the super harmonics. The first of these appeared in \cite{MR4258767}.

\begin{lemma}\label{lem:d_preserved}
  Let $d = \sum_{j=1}^n \partial_{g_j} \theta_j$ with $g_j \in K[x_1, \ldots, x_n]$ be a $G$-equivariant operator which strictly lowers $x$-degree. Then $d \colon \cSH_G \to \cSH_G$ preserves the super harmonics of $G$. In particular, we have a cochain complex
  \begin{equation}\label{eq:cochain_d}
    0 \to K \to \cSH_G^0 \too{d} \cSH_G^1 \too{d} \cdots \too{d} \cSH_G^n \to 0.
  \end{equation}
  
  \begin{proof}
    See \cite[Cor.~5.6]{MR4258767}. The argument in the next proof is very similar.
  \end{proof}

\end{lemma}

In particular, the generalized exterior derivatives $\dif_1, \ldots, \dif_r$ preserve the super harmonics $\cSH_G$. Moreover, we have the following ``dual'' result.

\begin{lemma}\label{lem:delta_preserved}
  Let $\delta = \sum_J g_J \partial_{\theta_J}$ be a $G$-equivariant operator which strictly lowers $\theta$-degree (i.e.~$g_\varnothing = 0$), where each $g_J \in K[x_1, \ldots, x_n]$ is linear. Then $\delta \colon \cSH_G \to \cSH_G$ preserves the super harmonics of $G$.
  
  \begin{proof}
    Suppose $\eta \in \cSH_G$, so $\partial_\omega \eta = 0$ for all $\omega \in \cJ_G$. We must show $\partial_\omega \delta\eta = 0$. Let $\omega = \sum_K h_K \theta_K$. We may suppose $\delta$ and $\omega$ are bi-homogeneous, so $|J|$ and $|K|$ are constant in the expansions, and that $\omega$ is $G$-invariant. Since $\partial_{\theta_J} \partial_{\theta_K} = (-1)^{|J||K|} \partial_{\theta_K} \partial_{\theta_J}$, we have
    \begin{align*}
      \partial_\omega \delta - (-1)^{|J||K|} \delta \partial_\omega
        &= \sum_{J, K} \partial_{h_K} \partial_{\theta_K} g_J \partial_{\theta_J} - (-1)^{|J||K|} g_J \partial_{\theta_J} \partial_{h_K} \partial_{\theta_K} \\
        &= \sum_{J, K} (\partial_{h_K} g_J - g_J \partial_{h_K}) \partial_{\theta_K} \partial_{\theta_J}.
    \end{align*}
    
    It is easily seen that $\partial_h g - g \partial_h \in K[\partial_{x_1}, \ldots, \partial_{x_n}]$ when $h, g \in K[x_1, \ldots, x_n]$ and $g$ is linear. Indeed, it may be reduced to the identity $\partial_{x_i}^a x_i - x_i \partial_{x_i}^a = a \partial_{x_i}^{a-1}$. Thus $\partial_\omega \delta - (-1)^{|J||K|} \delta \partial_\omega = \partial_\lambda$ for some $\lambda \in K[x_1, \ldots, x_n, \theta_1, \ldots, \theta_n]$. Since $\omega$ is $G$-invariant and $\delta$ is $G$-equivariant, $\lambda$ is $G$-invariant. Since $\delta$ strictly lowers $\theta$-degree, $\lambda$ has positive $\theta$-degree or is $0$, so $\lambda \in \cJ_G$, and we find $\partial_\omega \delta \eta = 0$.
  \end{proof}
\end{lemma}

In particular, $\dif^\dagger = \sum_{j=1}^n x_j \partial_{\theta_j}$ satisfies the hypotheses of \Cref{lem:delta_preserved} and hence preserves the harmonics $\cSH_G$. However, the adjoints $\dif_i^\dagger$ of the generalized exterior derivatives do not in general preserve the harmonics.

\subsection{Hodge theory and Laplacians}

A classic technique due to Hodge \cite[Ch.~III]{MR0003947} for analyzing cohomology on Riemannian manifolds involves replacing the cohomology groups with kernels of \textit{Laplacians}. An algebraic version of this decomposition for cell complexes was introduced by Eckmann \cite{MR13318} (see also \cite{MR1760592} and \cite[\S3]{MR1926878} for further references). We state the version of this decomposition appropriate to our context and include a standard proof sketch using elementary linear algebra for the benefit of the reader.

\begin{lemma}\label{lem:hodge}
  Suppose
  \begin{cd}
    A \rar[bend left]{\delta^\dagger}
      & B \rar[bend left]{d} \lar[bend left]{\delta}
      & C \lar[bend left]{d^\dagger}
  \end{cd}
  is a sequence of linear maps between finite-dimensional $K$-vector spaces with non-degenerate Hermitian forms, the adjoints are taken with respect to these forms, and $d \delta^\dagger = 0$. Then
  \begin{align}
    B &= \im\delta^\dagger \oplus \ker L \oplus \im d^\dagger \qquad\text{and}\label{lem:hodge.1}\\
    \ker d &= \im \delta^\dagger \oplus \ker L \label{lem:hodge.2}
  \end{align}
  where
    \[ L \coloneqq d^\dagger d + \delta^\dagger \delta \]
  is the ``total Laplacian''.
  
  \begin{proof}
    (Sketch.) Since $\ker d \perp \im d^\dagger$ and $B$ is finite-dimensional, $B = \ker d \oplus \im d^\dagger$, so we must show \eqref{lem:hodge.2}. Write $a = d^\dagger d$, $\alpha = \delta^\dagger \delta$. Now \eqref{lem:hodge.2} becomes $\ker a = \im \alpha \oplus \ker(a+\alpha)$. Clearly $\ker a \cap \ker \alpha \subset \ker(a+\alpha)$. Since $a\alpha = \alpha a = 0$ and $a=a^\dagger, \alpha=\alpha^\dagger$, one finds $\ker(a+\alpha) = \ker a \cap \ker \alpha$. Hence we must show $\ker a = \im \alpha \oplus (\ker a \cap \ker \alpha)$.
    
    Since $\alpha = \alpha^\dagger$, we have $B = \im \alpha \oplus (\im \alpha)^\perp = \im \alpha \oplus \ker \alpha$, so $\ker a = (\im \alpha \oplus \ker \alpha) \cap \ker a$. Now $\im d^\dagger \subset \ker d$ implies $\im \alpha \subset \ker a$, so $(\im \alpha \oplus \ker \alpha) \cap \ker a = \im \alpha \oplus (\ker \alpha \cap \ker a)$.

  \end{proof}
\end{lemma}

\begin{corollary}\label{cor:hodge_homology}
  The homology of the sequence $A \too{\delta^\dagger} B \too{d} C$ from \Cref{lem:hodge} at $B$ is
    \[ H \coloneqq \frac{\ker d}{\im \delta^\dagger} \cong \ker L \]
  where $L = d^\dagger d + \delta^\dagger \delta$. In particular, the sequence is exact at $B$ if and only if $L$ is invertible.
\end{corollary}

When $d = \dif$ is the exterior derivative, the adjoint operator is $\dif^\dagger = \sum_{j=1}^n x_j \partial_{\theta_j}$, and
  \[ L = \dif^\dagger \dif + \dif \dif^\dagger = \sum_{i=1}^n \left(\theta_i \partial_{\theta_i} + x_i \partial_{x_i}\right). \]
Thus $L$ acts by multiplying by the total degree, so it is invertible except on constants. Many of the operators $\dif_i$ in \Cref{tab:explicit} are of the form $d = \sum_{j=1}^n \partial_{x_j^N} \theta_j$ for some $N \in \mathbb{Z}_{\geq 0}$, which remain ``largely'' invertible.

\begin{lemma}\label{lem:power_sum_Laplacian}
  Let $d = \sum_{j=1}^n \partial_{x_j^N} \theta_j$ for some $N \in \mathbb{Z}_{\geq 0}$. Then the total Laplacian $L = d^\dagger d + d d^\dagger \in \End_K(K[x_1, \ldots, x_n, \theta_1, \ldots, \theta_n])$ has
    \[ \ker L = \Span_K\{x^\alpha : \alpha_j < N \text{ for }j=1, \ldots, n\}. \]
  
  \begin{proof}
    We first show that $L$ acts diagonally on the monomial basis. Recall that $\theta_j \partial_{\theta_\ell} + \partial_{\theta_\ell} \theta_j = \delta_{j=\ell}$. We compute
    \begin{align*}
      L &= d^\dagger d + d d^\dagger
            = \sum_{j,\ell=1}^n x_\ell^N \partial_{\theta_\ell} \partial_{x_j^N} \theta_j + \partial_{x_j^N} \theta_j x_\ell^N \partial_{\theta_\ell} \\
         &= \sum_{j \neq \ell} x_\ell^N \partial_{x_j^N} (\partial_{\theta_\ell} \theta_j + \theta_j \partial_{\theta_\ell}) + \sum_{j=1}^n \left(x_j^N \partial_{\theta_j} \partial_{x_j^N} \theta_j + \partial_{x_j^N} \theta_j x_j^N \partial_{\theta_j}\right) \\
         &= 0 + \sum_{j=1}^n \left(x_j^N \partial_{x_j^N} \partial_{\theta_j} \theta_j + \partial_{x_j^N} x_j^N \theta_j \partial_{\theta_j}\right).
    \end{align*}
    Suppose that $I = \{i_1 < \cdots < i_k\}$. Write $\theta_I = \theta_{i_1} \cdots \theta_{i_k}$. We have $\partial_{\theta_j} \theta_j \theta_I = \delta_{j \not\in I} \theta_I$ and $\theta_j \partial_{\theta_j} \theta_I = \delta_{j \in I} \theta_I$. Thus
      \[ L x^\alpha \theta_I = \left(\sum_{j \not\in I} x_j^N \partial_{x_j^N} + \sum_{j \in I} \partial_{x_j^N} x_j^N\right) x^\alpha \theta_I. \]
    Write $c^{\underline{N}} \coloneqq c(c-1)\cdots(c-N+1)$ for the falling factorial. Now
      \[ L x^\alpha \theta_I = \left(\sum_{j \not\in I} \alpha_j^{\underline{N}} + \sum_{j \in I} (\alpha_j + N)^{\underline{N}}\right) x^\alpha \theta_I. \]
    If $I \neq \varnothing$, the coefficient is positive. If $I = \varnothing$, the coefficient is zero if and only if $\alpha_j < N$ for all $N$.
  \end{proof}
\end{lemma}

\subsection{Proving exactness and the rank $2$ case}

We may now restate and prove \Cref{thm:exact} from the introduction.

\begin{reptheorem}{thm:exact}
  For any pseudo-reflection group $G \subset \GL(V)$ with $r = \dim(V/V^G)$, the exterior derivative cochain complex
  \begin{equation*}
    0 \to K \to \cSH_G^0 \too{\dif} \cSH_G^1 \too{\dif} \cdots \too{\dif} \cSH_G^r \too{\dif} 0
  \end{equation*}
  is exact.

  \begin{proof}
    By \Cref{lem:d_preserved}, the exterior derivative $\dif$ preserves $\cSH_G$ and the complex is well-defined. By \Cref{lem:delta_preserved}, the adjoint $\dif^\dagger$ also preserves $\cSH_G$. Thus by \Cref{cor:hodge_homology}, exactness at $\cSH_G^k$ for $k>0$ is equivalent to invertibility of the total Laplacian $L = \dif^\dagger \dif + \dif \dif^\dagger$ acting on $\cSH_G^k$. In this case, $L$ is invertible by \Cref{lem:power_sum_Laplacian}, in either $K[x_1, \ldots, x_n, \theta_1, \ldots, \theta_n]$ or $\cSH_G^k$. Finally, $\ker \dif|_{\cSH_G^0}$ is the set of all harmonics $f \in K[x_1, \ldots, x_n]$ where $\partial_{x_j} f = 0$ for all $j$, which is precisely $K$.
  \end{proof}
\end{reptheorem}

\begin{corollary}\label{cor:dif_difdagger}
  For $i+k \geq 1$,
  \begin{equation}\label{eq:cor:dif_difdagger.1}
    \cSH_G^{i,k} = \dif \cSH_G^{i+1,k-1} \oplus \dif^\dagger \cSH_G^{i-1,k+1}.
  \end{equation}
  
  \begin{proof}
    Apply \Cref{thm:exact} and \Cref{lem:hodge}.
  \end{proof}
\end{corollary}

\begin{remark}
  For $i>1$, the adjoints $\dif_i^\dagger$ do not typically preserve $\cSH_G$. Hence we cannot simply use \Cref{lem:power_sum_Laplacian} to analyze the homology of the complex $(\cSH_G^\bullet, \dif_i)$. A possible approach to \Cref{conj:Hilb_alt} would be to quantify this failure.
\end{remark}

Before restating and proving \Cref{thm:A1}, we give a criterion for the second component in \eqref{eq:cor:dif_difdagger.1} to satisfy \eqref{eq:thm:A}. Here we use the notation
  \[ \cD_G \coloneqq K[\partial_{x_1}, \ldots, \partial_{x_n}] \cSH_G^{\det}. \]

\begin{lemma}\label{lem:difdagger_theta}
  Suppose $\cSH_G^k \subset \cD_G$. Then
    \[ \dif^\dagger \cSH_G^k \subset \cD_G \]
  if and only if
    \[ \partial_{\theta_i} (\cSH_G^k)^{\det} \subset \cD_G \]
  for all $i = 1, \ldots, n$.
  
  \begin{proof}
    First suppose that $\dif^\dagger \cSH_G^k \subset \cD_G$.  A straightforward computation yields $\partial_{x_i} \dif^\dagger - \dif^\dagger \partial_{x_i} = \partial_{\theta_i}$. Thus
    \begin{align*}
      \partial_{\theta_i} (\cSH_G^k)^{\det}
        &= (\partial_{x_i} \dif^\dagger - \dif^\dagger \partial_{x_i}) (\cSH_G^k)^{\det} \\
        &\in \partial_{x_i} \cD_G + \dif^\dagger \cSH_G^k \subset \cD_G.
    \end{align*}
    
    Now suppose $\partial_{\theta_i} (\cSH_G^k)^{\det} \subset \cD_G$ for all $i = 1, \ldots, n$. Let $\omega \in \cSH_G^k$. Since $\cSH_G^k \subset \cD_G$, it suffices to take $\omega = \partial_g \eta$ for some homogeneous $g \in K[x_1, \ldots, x_n]$ and $\eta \in (\cSH_G^k)^{\det}$. We show that $\dif^\dagger \omega \in \cD_G$ by induction on $\deg g$. In the base case when $g$ is constant, $\omega \in (\cSH_G^k)^{\det}$, so $\dif^\dagger \omega \in (\cSH_G^{k-1})^{\det} \subset \cD_G$. If $\deg g > 0$, write $g = \sum_{j=1}^n x_j g_j$. Then
    \begin{align*}
      d^\dagger \omega
        &= d^\dagger \partial_g \eta = \sum_{j=1}^n \dif^\dagger \partial_{x_j} \partial_{g_j} \eta = \sum_{j=1}^n (\partial_{x_j} \dif^\dagger - \partial_{\theta_j}) \partial_{g_j} \eta \\
        &= \sum_{j=1}^n \partial_{x_j} \dif^\dagger \partial_{g_j} \eta - \partial_{g_j} \partial_{\theta_j} \eta \in \sum_{j=1}^n \partial_{x_j} \cD_G + \partial_{g_j} \cD_G \subset \cD_G.
    \end{align*}
  \end{proof}
\end{lemma}

\begin{reptheorem}{thm:A1}
  Let $G \subset \GL(V)$ be a pseudo-reflection group with rank $r = \dim(V/V^G)$. Then if $r \leq 2$ and either $G = G(m, 1, n)$ or $G$ is real,
  \begin{equation*}
  \begin{split}
    \cSH_G
      &= K[\partial_{x_1}, \ldots, \partial_{x_n}] \cSH_G^{\det} \\
      &= \{\partial_g \dif_{i_1} \cdots \dif_{i_k} \Delta_G : g \in K[x_1, \ldots, x_n], i_j \in [r] \}.
  \end{split}
  \end{equation*}
  
  \begin{proof}
    The $\cSH_G^0$ component of $\cSH_G$ satisfies \eqref{eq:thm:A} by Steinberg's \Cref{thm:Steinberg} in the sense that $\cSH_G^0 \subset K[\partial_{x_1}, \ldots, \partial_{x_n}] \cSH_G^{\det} \eqqcolon \cD_G$. Hence the $\dif \cSH_G^0$ summand of $\cSH_G^1$ from \Cref{cor:dif_difdagger} satisfies \eqref{eq:thm:A} as well. The result follows for $r \leq 1$, so take $r=2$.
    
    First suppose $G$ is real. By \Cref{thm:A2}, $\cSH_G^2 = \Span_K\{\vol\}$ where $\vol$ is the volume form on $V/V^G$, which transforms by $\det$. Since $\dif^\dagger$ is $G$-equivariant, $\dif^\dagger\cSH_G^2 \subset (\cSH_G^1)^{\det} \subset \cD_G$. Hence $\cSH_G^1 \subset \cD_G$ by \Cref{cor:dif_difdagger}, and the result follows.
    
    Now suppose $G = G(m, 1, n)$ has rank $r=2$, so $n=2$ and $m \geq 3$. By \Cref{thm:A2}, $\cSH_G^2 \subset \cD_G$. By \Cref{lem:difdagger_theta} and \Cref{thm:semi-invariants}, we may check that $\dif^\dagger\cSH_G^2 \subset \cD_G$ by instead checking that $\partial_{\theta_i} \dif_1 \dif_2 \Delta_G \in \cD_G$. We verify this condition directly for $G = G(m, 1, 2)$.
    
    From \Cref{tab:explicit}, we have
    \begin{align*}
      \Delta_G &= x_1^{m-1} x_2^{2m-1} - x_1^{2m-1} x_2^{m-1} \\
      \dif_1 &= \partial_{x_1} \theta_1 + \partial_{x_2} \theta_2 \\
      \dif_2 &= \partial_{x_1^{m+1}} \theta_1 + \partial_{x_2^{m+1}} \theta_2.
    \end{align*}
    Hence $(\cSH_G^2)^{\det} = K\,x_1^{m-2} x_2^{m-2} \theta_1 \theta_2$. Thus we must show $x_1^{m-2} x_2^{m-2} \theta_i \in \cD_G$. We directly compute
      \[ \partial_{x_1} (\dif_2 + \partial_{x_2^m} \dif_1) \Delta_G = 2 (2m-1)^{\underline{m+1}} (m-1) x_1^{m-2} x_2^{m-2} \theta_2. \]
    We may obtain $x_1^{m-2} x_2^{m-2} \theta_1$ symmetrically, which completes the proof.
  \end{proof}
\end{reptheorem}

\section{Gr\"obner and Artin bases for $G(m, p, n)$}\label{sec:artin}

Our proofs of \Cref{thm:B} and \Cref{thm:C} will use explicit Gr\"{o}bner and monomial bases of the coinvariant algebras $\cR_{G(m, p, n)}$ with respect to the lexicographic order on monomials. Artin \cite[p.41]{MR1616156} implicitly gave the first monomial basis in type $A$, which corresponds in a natural way to the \textit{inversion} statistic on permutations. Garsia \cite{MR597728} gave a separate \textit{descent basis} in type $A$ which corresponds to the \textit{major index} statistic. The descent basis was subsequently generalized to Weyl groups and $G(m, p, n)$ by a variety of authors; see \cite[p.324]{MR2342500} for details and further references. 

The Artin bases are very well-known for $\fS_n$, they appear to be folklore for $G(m, 1, n)$, and we have been unable to locate a description of them for $G(m, p, n)$ with $p>1$. We give Artin bases $\cA(m, p, n)$ for general $G(m, p, n)$ in \Cref{ssec:intro_artin_grobner}; see \Cref{def:Artin_mpn}. In \Cref{ssec:artin}, we give some combinatorial properties of $\cA(m, p, n)$. In \Cref{ssec:grobner}, we use Gr\"obner bases to prove the main results of this section, \Cref{thm:Artin_mpn} and  \Cref{thm:grobner_mpn}. Our arguments are elementary and self-contained, except for an appeal to Chevalley's result that $\dim_K \cR_G = |G|$ when $G$ is a (pseudo-)reflection group.

\subsection{Artin and Gr\"{o}bner bases of $\cR_{G(m, p, n)}$}\label{ssec:intro_artin_grobner}

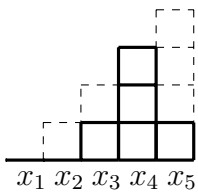
\begin{figure}[ht]
  \centering
  \begin{tikzpicture}[scale=0.5]
    \draw[dashed] (0, 0) grid (1, 0);
    \draw[dashed] (1, 0) grid (2, 1);
    \draw[dashed] (2, 0) grid (3, 2);
    \draw[dashed] (3, 0) grid (4, 3);
    \draw[dashed] (4, 0) grid (5, 4);
    
    \draw[very thick] (0, 0) grid (1, 0);
    \draw[very thick] (1, 0) grid (2, 0);
    \draw[very thick] (2, 0) grid (3, 1);
    \draw[very thick] (3, 0) grid (4, 3);
    \draw[very thick] (4, 0) grid (5, 1);
    
    \node[anchor=north west] at (0, 0) {$x_1$};
    \node[anchor=north west] at (1, 0) {$x_2$};
    \node[anchor=north west] at (2, 0) {$x_3$};
    \node[anchor=north west] at (3, 0) {$x_4$};
    \node[anchor=north west] at (4, 0) {$x_5$};
  \end{tikzpicture}
  \caption{A vertically oriented sub-staircase diagram visualizing the monomial $x_1^0 x_2^0 x_3^1 x_4^3 x_5^1$ in the Artin basis for $\fS_5$.}
  \label{fig:substaircase.1}
\end{figure}

The classical Artin basis for $\cR_n = \bC[x_1, \ldots, x_n]/\langle e_1, \ldots, e_n\rangle$ comes from \cite[p.41]{MR1616156} and consists of the $n!$ monomials
\begin{equation}
  \cA(1, 1, n) \coloneqq \{x_1^{a_1} \cdots x_n^{a_n} : 0 \leq a_i < i, \forall i \in [n]\}.
\end{equation}
These monomials may be visualized as ``sub-staircase diagrams''; see \Cref{fig:substaircase.1}.

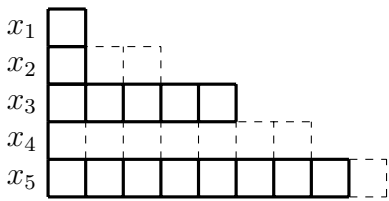
\begin{figure}[ht]
  \centering
  \begin{tikzpicture}[scale=0.5]
    \draw[dashed] (0, 0) grid (1, -1);
    \draw[dashed] (0, -1) grid (3, -2);
    \draw[dashed] (0, -2) grid (5, -3);
    \draw[dashed] (0, -3) grid (7, -4);
    \draw[dashed] (0, -4) grid (9, -5);

    \draw[very thick] (0, 0) grid (1, -1);
    \draw[very thick] (0, -1) grid (1, -2);
    \draw[very thick] (0, -2) grid (5, -3);
    \draw[very thick] (0, -3) grid (0, -4);
    \draw[very thick] (0, -4) grid (8, -5);
    
    \node[anchor=north east] at (0, 0) {$x_1$};
    \node[anchor=north east] at (0, -1) {$x_2$};
    \node[anchor=north east] at (0, -2) {$x_3$};
    \node[anchor=north east] at (0, -3) {$x_4$};
    \node[anchor=north east] at (0, -4) {$x_5$};
  \end{tikzpicture}
  \caption{A horizontally oriented sub-staircase diagram visualizing the monomial $x_1^{a_1} \cdots x_5^{a_5} = x_1 x_2 x_3^5 x_5^8$ in the Artin basis for $\fB_5 = G(2, 1, 5)$.}
  \label{fig:substaircase.2}
\end{figure}

The corresponding Artin basis for $G(m, 1, n)$ consists of the image of the $m^n n!$ monomials
\begin{equation}\label{eq:artin_mn}
  \cA(m, 1, n) \coloneqq \{x_1^{a_1} \cdots x_n^{a_n} : 0 \leq a_i < im, \forall i \in [n]\}.
\end{equation}
These too may be visualized as sub-staircase diagrams; see \Cref{fig:substaircase.2}. We typically draw these sub-staircase diagrams horizontally rather than vertically for convenience.

\begin{definition}
  A \textit{sub-staircase diagram of type $G(m, 1, n)$} is a left-justified arrangement of $n$ rows consisting of $a_1, \ldots, a_n$ square cells from top to bottom satisfying $0 \leq a_i < im$ for all $i \in [n]$.
\end{definition}

The Artin basis for general $G(m, p, n)$ is an index $p$ subset of $\cA(m, 1, n)$ which we may describe using the following operation on sub-staircase diagrams.

\begin{definition}\label{def:p-contraction}
  Suppose $m, p, n \in \bZ_{\geq 1}$ with $p \mid m$. Let $A$ be a sub-staircase diagram for $G(m, 1, n)$. Define the \textit{$p$-contraction} of $A$ as follows. Let $i$ be the largest index such that the $i$th row of $A$ has fewer than $m$ cells. Take the lower-left rectangle of width $m$ using rows $i, i+1, \ldots, n$ and shrink this rectangle horizontally by a factor of $p$. If this creates a partial cell, delete it. The result is the $p$-contraction of $A$.
\end{definition}

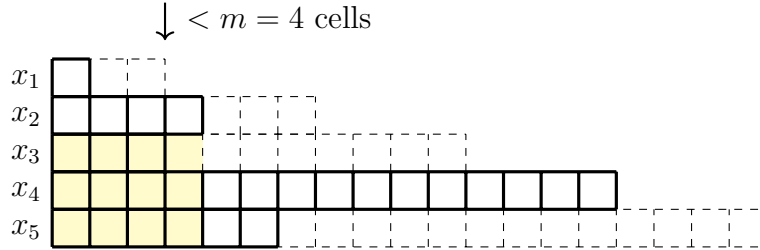
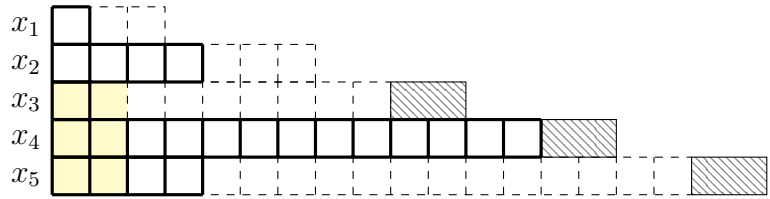
\begin{figure}[ht]
  \centering
  \begin{subfigure}[t]{\textwidth}
    \centering
    \begin{tikzpicture}[scale=0.5]
      \fill[yellow!25] (0, -2) rectangle (4, -5);

      \draw[dashed] (0, 0) grid (3, -1);
      \draw[dashed] (0, -1) grid (7, -2);
      \draw[dashed] (0, -2) grid (11, -3);
      \draw[dashed] (0, -3) grid (15, -4);
      \draw[dashed] (0, -4) grid (19, -5);

      \draw[very thick] (0, 0) grid (1, -1);
      \draw[very thick] (0, -1) grid (4, -2);
      \draw[very thick] (0, -2) grid (3, -3);
      \draw[very thick] (0, -3) grid (15, -4);
      \draw[very thick] (0, -4) grid (6, -5);
    
      \node[anchor=north east] at (0, 0) {$x_1$};
      \node[anchor=north east] at (0, -1) {$x_2$};
      \node[anchor=north east] at (0, -2) {$x_3$};
      \node[anchor=north east] at (0, -3) {$x_4$};
      \node[anchor=north east] at (0, -4) {$x_5$};
      
      \draw[->, thick] (3, 1.5) -- (3, 0.5) node[anchor=south west]{\ $< m = 4$ cells};
    \end{tikzpicture}
    \caption{A sub-staircase diagram for $G(4, 1, 5)$ representing the monomial $x_1 x_2^4 x_3^3 x_4^{15} x_5^6$. The bottom-most row with $< m$ cells is row $i=3$. The lower-left rectangle involved in $p$-contraction is highlighted.}
  \end{subfigure}
  \vspace{1em}
  
  \begin{subfigure}[t]{\textwidth}
    \centering
    \begin{tikzpicture}[scale=0.5]
      \fill[yellow!25] (0, -2) rectangle (2, -5);

      \draw[dashed] (0, 0) grid (3, -1);
      \draw[dashed] (0, -1) grid (7, -2);
      \draw[dashed] (0, -2) grid (9, -3);
      \draw[dashed] (0, -3) grid (13, -4);
      \draw[dashed] (0, -4) grid (17, -5);

      \draw[pattern=north west lines, pattern color=gray] (9, -2) rectangle (11, -3);
      \draw[pattern=north west lines, pattern color=gray] (13, -3) rectangle (15, -4);
      \draw[pattern=north west lines, pattern color=gray] (17, -4) rectangle (19, -5);

      \filldraw[very thick, fill=gray!25] (0, 0) grid (1, -1);
      \draw[very thick] (0, -1) grid (4, -2);
      \draw[very thick] (0, -2) grid (1, -3);
      \draw[very thick] (0, -3) grid (13, -4);
      \draw[very thick] (0, -4) grid (4, -5);
    
      \node[anchor=north east] at (0, 0) {$x_1$};
      \node[anchor=north east] at (0, -1) {$x_2$};
      \node[anchor=north east] at (0, -2) {$x_3$};
      \node[anchor=north east] at (0, -3) {$x_4$};
      \node[anchor=north east] at (0, -4) {$x_5$};
    \end{tikzpicture}
    \caption{The result of $2$-contracting the previous sub-staircase diagram. The half-cell arising from contracting row $3$'s three cells by a factor of two has been removed to create a sub-staircase diagram of type $G(4, 2, 5)$ representing the monomial $x_1 x_2^4 x_3 x_4^{13} x_5^4$.}
  \end{subfigure}
  \caption{An example of $p$-contraction.}
  \label{fig:p-contraction}
\end{figure}

\noindent Note that such an $i$ exists since the first row has length $< m$. See \Cref{fig:p-contraction} below for an example.

\begin{definition}\label{def:Artin_mpn}
  A \textit{sub-staircase diagram} of type $G(m, p, n)$ is the result of applying $p$-contraction to any sub-staircase diagram of type $G(m, 1, n)$. The \textit{Artin basis} of $G(m, p, n)$ is the corresponding set of monomials:
  \begin{equation}\label{eq:Artin_mpn}
      \begin{split}
        \cA(m, p, n) \coloneqq
        \{x_1^{a_1} \cdots x_n^{a_n} : \exists j \in [n] \text{ s.t. }
          &0 \leq a_i < im, \forall i < j, \text{ and } \\
          &0 \leq a_j < \frac{m}{p}, \text{ and } \\
          &0 \leq a_i - \frac{m}{p} < (i-1)m, \forall i > j\}.
      \end{split}
  \end{equation}
\end{definition}

\begin{example}
  The dihedral group of order $2m$ is $G(m, m, 2)$. The $2m$ sub-staircase diagrams of type $G(m, m, 2) = \fDih_{2m}$ are
    \[ \{(0, 0), \ldots, (m-1, 0), (0, 1), \ldots, (0, m)\} \]
  with corresponding Artin basis
    \[ \cA(m, m, 2) = \{1, x, \ldots, x^{m-1}, y, \ldots, y^m\} \] where
  $x \coloneqq x_1$ and $y \coloneqq x_2$.
\end{example}

In \Cref{ssec:grobner} we will show that $\cA(m, p, n)$ is indeed a monomial basis:

\begin{theorem}\label{thm:Artin_mpn}
  The Artin basis $\cA(m, p, n)$ descends to a basis for the coinvariant algebra $\cR_{G(m, p, n)} = K[x_1, \ldots, x_n]/\cI_{G(m, p, n)}$.
\end{theorem}

The Gr\"obner basis for the type $A$ coinvariant ideal relative to lexicographic order (actually relative to any linear order with$x_1 > \cdots > x_n$) is well-known and can be extracted from Artin's original argument in \cite{MR1616156} with some effort. Let $h_j(x_1, \ldots, x_n)$ denote the \textit{complete homogeneous symmetric polynomial} of degree $j$ in $n$ variables. The general Gr\"obner basis is as follows, which is proved in \Cref{ssec:grobner}.

\begin{theorem}\label{thm:grobner_mpn}
  The reduced Gr\"obner basis of $\cI_{G(m, p, n)}$ with respect to the lexicographic term order with $x_1 > \cdots > x_n$ is:
  \begin{itemize}
    \item If $p=1$: $\{h_j(x_j^m, \ldots, x_n^m) : j \in [n]\}$.
    \item If $p>1$:
      \begin{align*}
        \{h_j(x_j^m, &\ldots, x_n^m) : j \in [n-1]\} \\
        &\sqcup \{h_{j-1}(x_j^m, \ldots, x_n^m) (x_j \cdots x_n)^{m/p} : j \in [n]\}.
      \end{align*}
  \end{itemize}
\end{theorem}

\begin{remark}
  Recall that a Gr\"obner basis $\cG \subset F[x_1, \ldots, x_n]$ is \textit{reduced} if for all $g \in \cG$, the leading coefficient of $g$ is $1$ and for all monomials $m$ in all $h \in \cG - \{g\}$, the leading monomial of $g$ does not divide $m$. See \cite[\S2.4, 2.7]{MR3330490} for details.
\end{remark}

\subsection{Combinatorial properties of the Artin bases for $G(m, p, n)$}\label{ssec:artin}

We now give some combinatorial properties of the Artin bases $\cA(m, p, n)$.

\begin{lemma}\label{lem:Artin_mpn_Hilb}
  There are $m^n n!/p$ sub-staircase diagrams of type $G(m, p, n)$. Moreover, the Hilbert series of the Artin basis of type $G(m, p, n)$ is
  \begin{equation}\label{eq:Artin_mpn_Hilb}
    \Hilb(\cA(m, p, n); q) = [m]_q [2m]_q \cdots [(n-1)m]_q [nm/p]_q.
  \end{equation}
  
  \begin{proof}
    Note that $p$-contraction results in some fraction $k/p$ of a cell, $k \in \{0, \ldots, p-1\}$, which is then discarded. It follows that the fibers of the $p$-contraction map are all of cardinality $p$. Hence there are $m^n n!/p$ sub-staircase diagrams of type $G(m, p, n)$.

    For the Hilbert series, condition on which row $j$ is the lowest with length $< m$. The resulting Hilbert series is easily seen to be
    \begin{align*}
      \sum_{j=1}^m &[m]_q [2m]_q \cdots [(j-1)m]_q [m/p]_q
        q^{m/p} [jm]_q \cdots q^{m/p} [(n-1)m]_q \\
        &= \sum_{j=1}^n \left(\prod_{i=1}^{n-1} [im]_q\right) q^{(n-j)m/p} [m/p]_q \\
        &= \prod_{i=1}^{n-1} [im]_q \cdot [m/p]_q \cdot \sum_{j=1}^n q^{(n-j) m/p}.
    \end{align*}
    The result follows by observing
    \begin{align*}
      [m/p]_q \cdot \sum_{j=1}^n q^{(n-j) m/p}
        &= \frac{1 - q^{m/p}}{1 - q} \cdot \frac{1 - q^{nm/p}}{1 - q^{m/p}} \\
        &= [nm/p]_q.
    \end{align*}
  \end{proof}
\end{lemma}

\begin{lemma}\label{lem:hooks}
  A sub-staircase diagram of type $G(m, 1, n)$ is a sub-staircase diagram of type
  $G(m, p, n)$ if and only if it does not contain any of the ``hook'' diagrams
  
  
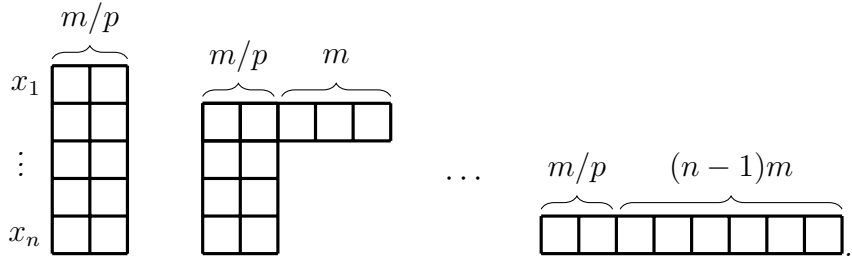
\begin{tikzpicture}[scale=0.5]
    \draw [decorate,decoration={brace,amplitude=5pt,raise=0.5ex}] (0,0) -- (2,0) node[midway,yshift=1.5em]{$m/p$};
    \draw[very thick] (0, 0) grid (2, -5);
    
    \node[anchor=north east] at (0, 0) {$x_1$};
    \node[anchor=north east, xshift=-0.5em, yshift=0.5em] at (0, -2) {$\vdots$};
    \node[anchor=north east] at (0, -4) {$x_n$};
    
    \begin{scope}[shift={(4, 0)}]
      \draw [decorate,decoration={brace,amplitude=5pt,raise=0.5ex}] (0,-1) -- (1.9,-1) node[midway,yshift=1.5em]{$m/p$};
      \draw [decorate,decoration={brace,amplitude=5pt,raise=0.5ex}] (2.1,-1) -- (5,-1) node[midway,yshift=1.5em]{$m$};
      \draw[very thick] (0, -1) grid (2, -5);
      \draw[very thick] (0, -1) grid (5, -2);
    \end{scope}
    
    \begin{scope}[shift={(11, 0)}]
      \node[very thick] at (0, -3) {$\cdots$};
    \end{scope}

    \begin{scope}[shift={(13, 0)}]
      \draw [decorate,decoration={brace,amplitude=5pt,raise=0.5ex}] (0,-4) -- (1.9,-4) node[midway,yshift=1.5em]{$m/p$};
      \draw [decorate,decoration={brace,amplitude=5pt,raise=0.5ex}] (2.1,-4) -- (8,-4) node[midway,yshift=1.5em]{$(n-1)m$};
      \draw[very thick] (0, -4) grid (8, -5);
    \end{scope}
  \end{tikzpicture}.
  
  \begin{proof}
    First suppose to the contrary that there exists a sub-staircase diagram $A$ of type $G(m, p, n)$ containing the $j$th hook diagram. We may reverse the $p$-contraction process by horizontally expanding the width $m/p$ rectangle in $A$ involving rows $i, i+1, \ldots, n$ where $i$ is the lowest row of length $<m/p$. The result, call it $B$, necessarily expands the rectangle of width $m/p$ involving rows $j, j+1, \ldots, n$ in the hook diagram. After expansion, row $j$ of $B$ has length $\geq m + (j-1)m = jm$, contradicting the fact that $B$ must be a sub-staircase diagram of type $G(m, 1, n)$. Thus sub-staircase diagrams of type $G(m, p, n)$ do not contain the above hooks.

    Conversely, suppose we have a sub-staircase diagram $A = (a_1, \ldots, a_n)$ of type $G(m, 1, n)$ which does not contain the above hooks. Since $A$ does not contain the first hook, an $m/p$ by $n$ rectangle, there is some row in $A$ of length $< m/p$. Thus we may apply the reverse $p$-contraction process described above to $A$ starting at some row $j$. We must only show the result, call it $B = (b_1, \ldots, b_n)$, remains a sub-staircase diagram of type $G(m, 1, n)$, i.e.~$b_i < im$ for all $i \in [n]$. We have $b_1 = a_1, \ldots, b_{j-1} = a_{j-1}$, so the necessary constraint is satisfied for $i < j$. When $i=j$, we have $b_j = a_jp < (m/p)p = m \leq jm$ as required. Finally, since rows $i > j$ of $A$ have length $\geq m/p$ and yet the $i$th hook is not contained in the diagram, we must have $a_i < (m/p) + (i-1)m$ for $i > j$. Thus
      \[ b_i = a_i + (m - m/p) < (i-1)m + m = im, \]
    as required.
  \end{proof}
\end{lemma}

We also have the following recursive description of the Artin bases.
The Hilbert series formula \eqref{eq:Artin_mpn_Hilb} also follows easily from
this description and induction.

\begin{lemma}
  The Artin bases of types $G(m, p, n)$ are defined recursively by
  \begin{align*}
    \cA(m, p, n) &= \{x_1^{a_1} \cdots x_n^{a_n} :
      \forall i \in [n-1], 0 \leq a_i < im \text{ and } 0 \leq a_n < m/p \} \sqcup \\
      &\bigsqcup_{j=m/p}^{(n-1)m + m/p - 1} x_n^j \cA(m, p, n-1)
  \end{align*}
  for $n \geq 2$, with base cases
    \[ \cA(m, p, 1) = \{x_1^{a_1} : 0 \leq a_1 < m/p\}. \]
  
  \begin{proof}
    The base cases are immediate. For the recursive formula,
    the first term consists of elements where $j=n$ when performing
    $p$-contraction. When $j<n$, we may keep or remove the final
    row without affecting the procedure materially, and that row has length satisfying
    $0 \leq a_n-m/p < (n-1)m$.
  \end{proof}
\end{lemma}

\subsection{Gr\"obner bases for $G(m, p, n)$}\label{ssec:grobner}

We now prove \Cref{thm:Artin_mpn} and \Cref{thm:grobner_mpn}.

\begin{lemma}\label{lem:grobner_mpn}
  The polynomials in the proposed Gr\"obner bases in \Cref{thm:grobner_mpn} belong to $\cI_{m, p, n}$.
  
  \begin{proof}
    We first prove by induction on $j$ that $h_k(x_j^m, \ldots, x_n^m) \in \cI_{m, p, n}$ whenever $k \geq j$. In the base case $j=1$, $h_k(x_1^m, \ldots, x_n^m)$ is $G(m, p, n)$-invariant, so $h_k(x_1^m, \ldots, x_n^m) \in \cI_{m, p, n}$ for all $k \geq 1$.
    
    For the inductive step, first recall
    \begin{equation}\label{eq:hk}
      h_{k+1}(x_j, \ldots, x_n) = h_{k+1}(x_{j+1}, \ldots, x_n) + x_j h_k(x_j, \ldots, x_n).
    \end{equation}
    Suppose $h_k(x_j, \ldots, x_n) \in \cI_{m, p, n}$ for all $k \geq j$. Since $h_{k+1}(x_j, \ldots, x_n)$ and $h_k(x_j, \ldots, x_n)$ belong to $\cI_{m, p, n}$, \eqref{eq:hk} gives $h_{k+1}(x_{j+1}, \ldots, x_n) \in \cI_{m, p, n}$, completing the induction. If $p=1$, we are now done, so take $p>1$.
    
    We next prove by induction on $j$ that $h_{j-1}(x_j^m, \ldots, x_n^m)(x_j \cdots x_n)^{m/p} \in \cI_{m, p, n}$. In the base case $j=1$, it is easy to see directly that $(x_1 \cdots x_n)^{m/p}$ is $G(m, n, p)$-invariant. For the inductive step, suppose
      \[ h_{j-1}(x_j^m, \ldots, x_n^m)(x_j \cdots x_n)^{m/p} \in \cI_{m, p, n}. \]
    Letting $k=j-1$ and $x_i \mapsto x_i^m$ in \eqref{eq:hk}, multiplying by $(x_{j+1} \cdots x_n)^{m/p}$, and rearranging gives
    \begin{align*}
      h_j(x_{j+1}^m, \ldots, x_n^m) &(x_{j+1} \cdots x_n)^{m/p}
        = (x_{j+1} \cdots x_n)^{m/p} h_j(x_j^m, \ldots, x_n^m)\\
        &- x_j^{m - m/p} (x_j \cdots x_n)^{m/p} h_{j-1} (x_j^m, \ldots, x_n^m)
        \in \cI_{m, n, p}.
    \end{align*}
  \end{proof}   
\end{lemma}

\begin{lemma}\label{lem:Artin_LT}
  The leading monomials of the polynomials in the proposed Gr\"obner basis in \Cref{thm:grobner_mpn} with respect to the lexicographic term order with $x_1 > \cdots > x_n$ are as follows.
  \begin{itemize}
    \item If $p=1$: $\{x_j^{jm} : j \in [n]\}$.
    \item If $p>1$: $\{x_j^{jm} : j \in [n-1]\}
      \sqcup \{x_j^{(j-1)m} (x_j \cdots x_n)^{m/p} : j \in [n]\}$.
  \end{itemize}

  \begin{proof}
    Clear.
  \end{proof}
\end{lemma}

We may now complete the proof of \Cref{thm:Artin_mpn} and \Cref{thm:grobner_mpn}.

\begin{proof}[{Proof of \Cref{thm:Artin_mpn}}]
  Let $\LT(S)$ denote the ideal generated by leading terms of $S \subset K[x_1, \ldots, x_n]$. Consider the monomial ideal $\LT(\cG)$ generated by the leading terms of the proposed Gr\"obner basis $\cG$ for $G(m, p, n)$. By \Cref{lem:grobner_mpn}, $\LT(\cG) \subset \LT(\cI_{m, p, n})$, so
  \begin{equation}\label{eq:grobner_mpn.1}
    K[x_1, \ldots, x_n]/\LT(\cG) \supset K[x_1, \ldots, x_n]/\LT(\cI_{m, p, n}).
  \end{equation}

  The standard monomial basis of $K[x_1, \ldots, x_n]/\LT(\cG)$ consists of all monomials not divisible by any of the leading terms of $\cG$. Encoding monomials in diagrams as in \Cref{ssec:intro_artin_grobner}, \Cref{lem:hooks} and \Cref{lem:Artin_LT} together show that the Artin basis $\cA(m, p, n)$ is precisely the standard monomial basis for $K[x_1, \ldots, x_n]/\LT(\cG)$. By \Cref{lem:Artin_mpn_Hilb}, this basis has size $m^n n!/p$. 

  As for the right-hand side of \eqref{eq:grobner_mpn.1}, a set of representatives of a basis of $K[x_1, \ldots, x_n]/\LT(\cI_{m, p, n})$ descends to a basis of $K[x_1, \ldots, x_n]/\cI_{m, p, n}$ as usual, and in particular their dimensions agree. Chevalley \cite{MR72877} proved
    \[ \dim K[x_1, \ldots, x_n]/\cI_{m, p, n} = |G(m, p, n)| = m^n n!/p. \]
  Thus equality must hold in \eqref{eq:grobner_mpn.1}, so $\LT(\cG) = \LT(\cI_{m, p, n})$ and the result follows.
\end{proof}

\begin{proof}[{Proof of  \Cref{thm:grobner_mpn}}]
  By the preceding argument, $\LT(\cG) = \LT(\cI_{m, p, n})$, so the proposed Gr\"obner basis is in fact a Gr\"obner basis. All that is left is to check that $\cG$ is reduced. The leading coefficients are all $1$, so we must only verify the divisibility condition.
  
  We begin with the $p=1$ case. Consider a general term
    \[ x_j^{a_j m} \cdots x_n^{a_n m} \text{ in }h_j(x_j^m, \ldots, x_n^m) \]
  where $j \in [n]$. Suppose to the contrary $x_k^{km} \mid x_j^{a_j m} \cdots x_n^{a_n m}$ for some $k \in [n]$ with $k \neq j$. Then $a_k m \geq km$, so $a_k \geq 1$, and hence $j \leq k$. Now $h_j(x_j^m, \ldots, x_n^m)$ has degree $jm$ while $x_k^{km}$ has degree $km$, so $km \leq jm$, forcing $k=j$, a contradiction.
  
  We now turn to the case $p>1$. Again consider a general term $x_j^{a_j m} \cdots x_n^{a_n m}$ in $h_j(x_j^m, \ldots, x_n^m)$ where $j \in [n-1]$. By the above argument, $x_k^{km} \nmid a_j^{a_j m} \cdots x_n^{a_n m}$ for all $k \in [n-1]$. Now suppose to the contrary that $x_k^{(k-1)m} (x_k \cdots x_n)^{m/p} \mid x_j^{a_j m} \cdots x_n^{a_n m}$ for some $k \in [n]$. Since $a_j m + \cdots + a_n m = jm$, we have $a_k \leq j$. As before, $j \leq k$. On the other hand, $(k-1)m + m/p \leq a_k m$, so $a_k > k-1$, hence $a_k \geq k \geq j$. Thus we have $a_k=k=j$ and $x_j^{a_j m} \cdots x_n^{a_n m} = x_k^{k m}$. However, $k < n$, so $x_k^{(k-1)m} (x_k \cdots x_n)^{m/p} \nmid x_k^{k m}$, as required.
  
  Finally, consider a general term
    \[ x_j^{a_j m} \cdots x_n^{a_n m} (x_j \cdots x_n)^{m/p} \text{ in } h_{j-1}(x_j^m, \ldots, x_n^m) (x_j \cdots x_n)^{m/p} \]
  for $j \in [n]$. Note that $a_j m + \cdots + a_n m = (j-1) m$, so $a_k \leq j-1$. First, suppose to the contrary that $x_k^{k m} \mid x_j^{a_j m} \cdots x_n^{a_n m} (x_j \cdots x_n)^{m/p}$ for $k \in [n-1]$. Then $km \leq a_k m + m/p$, so $k \leq a_k$ since $p>1$. Since $a_k \geq 1$, we have $j \leq k$. Now $j \leq k \leq a_k \leq j-1$, a contradiction. Finally, suppose to the contrary that $x_k^{(k-1) m} (x_k \cdots x_n)^{m/p} \mid x_j^{a_j m} \cdots x_n^{a_n m} (x_j \cdots x_n)^{m/p}$ for some $k \in [n]$ with $k \neq j$. Clearly $j \leq k$. Also, $(k-1)m + m/p \leq a_k m + m/p$, so $k-1 \leq a_k \leq j-1$, meaning $k \leq j$. Hence $j=k$, a final contradiction.
\end{proof}

%\subsection{Generating sets for $\cI_{G(m, p, n)}$}
%
%The following generating sets for the ideals $\cI_{m, p, n}$ are very well known, especially when $p=1$. We include a proof since it is very straightforward given the preceding material. Let
%  \[ p_j(x_1, \ldots, x_n) \coloneqq x_1^j + \cdots + x_n^j \]
%be the usual degree $j$ \textit{power-sum symmetric polynomial}.
%
%\begin{lemma}\label{lem:p_gens}
%  The ideal $\cI_{m, p, n}$ is generated by the following set.
%  \begin{itemize}
%    \item If $p=1$: $\{p_j(x_1^m, \ldots, x_n^m) : j \in [n]\}$.
%    \item If $p>1$: $ \{p_j(x_1^m, \ldots, x_n^m) : j \in [n-1]\} \sqcup
%      \{(x_1 \cdots x_n)^{m/p}\}$.
%  \end{itemize}
%  
%  \begin{proof}
%    Let $\cM$ denote the ideal generated by the proposed set. It is easy to check directly that each proposed generator is $G(m, p, n)$-invariant, so $\cM \subset \cI_{m, p, n}$. From Newton's identities, $\{p_j(x_1^m, \ldots, x_n^m) : j \in [N]\}$, $\{h_j(x_1^m, \ldots, x_n^m) : j \in [N]\}$, and $\{e_j(x_1^m, \ldots, x_n^m) : j \in [N]\}$ all generate the same subalgebra (over a base ring containing $\bQ$), where $e_j$ denotes an elementary symmetric polynomial and $N \leq n$. Since $(x_1 \cdots x_n)^m = e_n(x_1^m, \ldots, x_n^m)$, it follows that $\{h_j(x_1^m, \ldots, x_n^m) : j \in [n]\} \subset \cM$, even when $p > 1$. The arguments in \Cref{lem:grobner_mpn} now show that $\cM$ contains the Gr\"obner basis in \Cref{thm:grobner_mpn}. Thus $\cM = \cI_{m, p, n}$.
%  \end{proof}
%\end{lemma}

\section{Bi-degree bounds for $G(m, 1, n)$}\label{sec:bidegree}

We now prove one of our main results, \Cref{thm:B}. The argument uses the Gr\"obner bases in \Cref{thm:grobner_mpn}. We first record two recursive relations involving the elements $\partial_{x_i} h_j(x_1^m, \ldots, x_n^m)$.

\begin{lemma}\label{lem:recursive}
  Let $j, m \in \mathbb{Z}_{\geq 1}$. Then
  \begin{equation}\label{eq:recursive.rel.1}
    \begin{split}
    \partial_{x_i} h_j(x_1^m, \ldots, x_n^m)
      &= m x_i^{m-1} h_{j-1}(x_1^m, \ldots, x_n^m) \\
      &\quad+ x_i^m \partial_{x_i} h_{j-1}(x_1^m, \ldots, x_n^m).
    \end{split}
  \end{equation}
  and
  \begin{equation}\label{eq:recursive.rel.2}
    \begin{split}
    \sum_{i=1}^n &(x_1 \cdots \widehat{x}_i \cdots x_n)^{m-1} \partial_{x_i} h_j(x_1^m, \ldots, x_n^m) \\
      &= (x_1 \cdots x_n)^{m-1} m(n+j-1)h_{j-1}(x_1^m, \ldots, x_n^m).
    \end{split}
  \end{equation}
  
  \begin{proof}
    To prove \eqref{eq:recursive.rel.1}, start with
    \begin{equation}\label{eq:recursive.rel.3}
      h_j(x_1, \ldots, x_n) = h_j(x_1, \ldots, \widehat{x}_i, \ldots, x_n) + x_i h_{j-1}(x_1, \ldots, x_n),
    \end{equation}
    apply $x_\ell \mapsto x_\ell^m$ for all $1 \leq \ell \leq n$, and differentiate with respect to $x_i$. Equation \eqref{eq:recursive.rel.2} is implied by \eqref{eq:recursive.rel.1} and Euler's identity $\sum_{i=1}^n x_i \partial_{x_i} h = \deg(h) h$ for $h$ homogeneous.
  \end{proof}
\end{lemma}

The argument in the following key lemma is based on one given by Fran\c{c}ois Brunault \cite{MO98054} in his proof that the $h_j$ are typically irreducible.

\begin{lemma}\label{lem:Brunault}
  For each $j, m \in \bZ_{\geq 1}$, the $n$ polynomials
    \[ \{\partial_{x_i} h_j(x_1^m, \ldots, x_n^m) : i \in [n]\} \]
  have no common zero in $\bC^n - \{0\}$.
  
  \begin{proof}
    We argue first by induction on $n$, and for each $n$ by induction on $j$. The result is clear in the base cases $n=1$ or $j=1$, so take $n, j \geq 2$. Suppose $z = (z_1, \ldots, z_n)$ is a common zero. If $z_\ell = 0$ for some $\ell \in [n]$, then $(z_1, \ldots, \widehat{z}_\ell, \ldots, z_n)$ is a common zero of $\{\partial_{x_i} h_j(x_1^m, \ldots, \widehat{x}_\ell, \ldots, x_n^m) : i \in [n] - \{\ell\}\}$ by \eqref{eq:recursive.rel.3}, so $z=0$ by induction on $n$. Hence we suppose $z_\ell \neq 0$ for all $\ell \in [n]$.
    
    By \eqref{eq:recursive.rel.2} and the fact that $z_\ell \neq 0$, $h_{j-1}(z_1^m, \ldots, z_n^m) = 0$. Equation \eqref{eq:recursive.rel.1} implies that $(z_1, \ldots, z_n)$ is a zero of $\partial_{x_i} h_{j-1}(x_1^m, \ldots, x_n^m)$ for all $i \in [n]$. Thus $z=0$ by induction on $j$, completing the proof.
  \end{proof}
\end{lemma}

\begin{remark}
  \Cref{lem:Brunault} and the Nullstellensatz combine to show that the radical of $\langle\partial_{x_i} h_j(x_1^m, \ldots, x_n^m) : i \in [n]\rangle$ is $\langle x_1, \ldots, x_n\rangle$, i.e.
    \[ \{\partial_{x_i} h_j(x_1^m, \ldots, x_n^m)\} \]
  forms a \textit{homogeneous system of parameters} for $\mathbb{C}[x_1, \ldots, x_n]$ (for $j \geq 2$). Since $\bC[x_1, \ldots, x_n]$ is \textit{Cohen--Macaulay}, they form a \textit{regular sequence}. Consequently,
  \begin{equation}\label{eq:di_hk_Hilb}
    \Hilb\left(\frac{\bC[x_1, \ldots, x_n]}{\langle \partial_{x_i} h_j(x_1^m, \ldots, x_n^m) : i \in [n]\rangle}; q\right)
     = \left(\frac{1-q^{mj-1}}{1-q}\right)^n
     = [mj-1]_q^n,
  \end{equation}
  where the exponents in the numerator arise from the degrees of the elements in the regular sequence. See \cite[I.5, p.39-42]{MR725505} for details and \cite{MR2542141} for results similar to \Cref{lem:Brunault}.
\end{remark}

\begin{corollary}\label{cor:Brunault.2}
  Let $K \subset \mathbb{C}$. For each $j, m \in \bZ_{\geq 1}$ and homogeneous polynomial $p(x_j, \ldots, x_n) \in K[x_j, \ldots, x_n]$, we have
    \[ \deg p > (mj-2)(n-j+1) \quad\Rightarrow\quad
        p \in \langle \partial_{x_i} h_j(x_j^m, \ldots, x_n^m) : i=j, \ldots, n\rangle. \]
  
  \begin{proof}
    By \eqref{eq:di_hk_Hilb} and \Cref{lem:Brunault},
    \begin{align*}
      \Hilb\left(\frac{K[x_j, \ldots, x_n]}{\langle \partial_{x_i} h_j(x_j^m, \ldots, x_n^m)
        : i=j, \ldots, n\rangle}; q\right)
       &= \left(\frac{1-q^{mj-1}}{1-q}\right)^{n-j+1} \\
       &= [mj-1]_q^{n-j+1}.
    \end{align*}
    The right-hand side has degree $(mj-2)(n-j+1)$.
  \end{proof}
\end{corollary}

%\begin{remark}
%  We have in fact found an alternate, purely combinatorial approach to \Cref{cor:Brunault.2}, avoiding use of regular sequences and complex numbers. It would be interesting to extend either argument to the super diagonal coinvariant algebra.
%\end{remark}

The final ingredient needed for \Cref{thm:B} is the following observation. Afterwards we restate and prove \Cref{thm:B}.

\begin{lemma}\label{lem:pjhr_wedge}
  If $j \in [n]$ and $i \in [n]$, then
    \[ \partial_{x_i} h_j(x_j^m, \ldots, x_n^m)\,\theta_j \cdots \theta_n \in \cJ_{G(m, p, n)}. \]
  
  \begin{proof}
    The result is trivial if $i<j$, so take $i \geq j$. Since the exterior derivative $\dif = \sum_{\ell=1}^n \partial_{x_\ell} \theta_\ell$ is $G$-equivariant and an anti-derivation, it preserves $\cJ_{G(m, p, n)}$. By \Cref{thm:grobner_mpn},
      \[ h_j(x_j^m, \ldots, x_n^m) \in \cI_{G(m, p, n)} \subset \cJ_{G(m, p, n)}. \]
    Thus
    \begin{align*}
      \dif h_j(x_j^m, \ldots, &x_n^m)\,\theta_j \cdots \widehat{\theta}_i \cdots \theta_n \\
        &= \pm\partial_{x_i} h_j(x_j^m, \ldots, x_n^m)\,\theta_j \cdots \theta_n
          \in \cJ_{G(m, p, n)}.
    \end{align*}
  \end{proof}
\end{lemma}

\begin{reptheorem}{thm:B}
  Let $G = G(m, 1, n)$. Then $\cSH_{G(m, 1, n)}^{i, k} \neq 0$ if and only if $i, k \geq 0$ and
  \begin{equation*}
    i+k+m\binom{k}{2} \leq m \binom{n}{2} + (m-1)n.
  \end{equation*}

  \begin{proof}
    First suppose $i+k+m\binom{k}{2} > m \binom{n}{2} + (m-1)n$. Let $f \in \cSR_{G(m, 1, n)}^{i, k}$. We must show $f = 0$. It suffices to suppose $f = x^\alpha\,\theta_I$ is a monomial with $|\alpha| = i$ and $|I| = k$ and show $x^\alpha\,\theta_I \in \cJ_{G(m, 1, n)}$. Since $\cSR_{G(m, 1, n)}$ is $\fS_n$-invariant, we may take $\theta_I = \theta_j \cdots \theta_n$ where $j = n-k+1$. Since $\cJ_{G(m, 1, n)}$ contains $\cI_{G(m, 1, n)} \cdot K[\theta_1, \ldots, \theta_n]$, we may further suppose $x^\alpha$ belongs to the Artin basis $\cA(m, 1, n)$ for $\cR(m, 1, n)$, so $\alpha_\ell \leq \ell m - 1$ for all $\ell \in [n]$ by \Cref{thm:Artin_mpn}.
    
    The degree of $x_j^{\alpha_j} \cdots x_n^{\alpha_n}$ is
    \begin{align*}
      \alpha_j + \cdots + \alpha_n
        &= i - \alpha_1 - \cdots - \alpha_{j-1} \\
        &\geq i - \sum_{\ell=1}^{j-1} (\ell m-1)
          = i - m \binom{j}{2} + (j-1) \\
        &> m \binom{n}{2} + (m-1)n - k - m\binom{k}{2} - m\binom{j}{2} + (j-1) \\
        &= (mj-2)(n-j+1).
    \end{align*}
    By \Cref{cor:Brunault.2}, $x_j^{\alpha_j} \cdots x_n^{\alpha_n}$ belongs to the ideal generated by
      \[ \partial_{x_\ell} h_r(x_j^m, \ldots, h_n^m). \]
    Hence by \Cref{lem:pjhr_wedge},
      \[ x_j^{\alpha_j} \cdots x_n^{\alpha_n}\,\theta_j \cdots \theta_n
          \in \cJ_{G(m, 1, n)}. \]
    Thus $x^\alpha\,\theta_I \in \cJ_{G(m, 1, n)}$, which proves necessity.
    
    We now prove sufficiency. We first show that if $i+k+m\binom{k}{2} = m \binom{n}{2} + (m-1)n$, then $\cSH_{G(m, 1, n)}^{i,k} \neq 0$. By \Cref{tab:explicit}, $\deg \Delta_V = m\binom{n}{2} + (m-1)n$ and the operators $\dif_j$ lower $x$-degree by $m(j-1)+1$ while raising $\theta$-degree by $1$. By \Cref{thm:semi-invariants}, the element $\dif_1 \cdots \dif_k \Delta_V \in \cSH_G$ is a non-zero harmonic $k$-form of $x$-degree
    \begin{align*}
      m\binom{n}{2} + &(m-1)n - \sum_{j=1}^k (m(j-1) + 1) \\
        &= m\binom{n}{2} + (m-1)n - m\binom{k}{2} - k \\
        &= i.
    \end{align*}
    Thus $\cSH_{G(m, 1, n)}^{i,k} \neq 0$.
    
    Finally, since $\cSH_G$ is closed under partial differentiation, $\cSH_G^{i,k} \neq 0$ implies $\cSH_G^{i',k'} \neq 0$ for all $i' \leq i$ and $k' \leq k$. Indeed, one may pick a monomial $c x^\alpha \theta_I$ in a non-zero element $\omega \in \cSH_G^{i,k}$. Now $\partial_{x^\alpha \theta_I} \omega = \langle cx^\alpha \theta_I, x^\alpha \theta_I\rangle$ is a non-zero constant, so $0 \neq \partial_{x^\beta \theta_J} \omega \in \cSH_G^{i',k'}$ for $\beta \leq \alpha, J \subset K$ with $|\beta| = i'$, $|J| = k'$.
  \end{proof}
\end{reptheorem}

\section{Total degree bounds for $G(m, p, n)$}\label{sec:total}

We next prove \Cref{thm:C}. We begin with the case of $1$-forms.

\begin{lemma}\label{lem:k1_top}
  Let $G = G(m, p, n)$. The top-degree component of $\cSH_G^1$ is $K\,\dif \Delta_G$, except when $G = G(2, 2, 2)$ when it is
    \[ \Span_K \{\dif \Delta_G, \dif_2 \Delta_G\} = \Span_K\{x_1 \theta_1 + x_2 \theta_2, x_2 \theta_1 + x_1 \theta_2\}. \]
  
  \begin{proof}
    Suppose $\omega \in \cSH_G^1$ is homogeneous with $\omega = \sum_{j=1}^n \omega_j \theta_j$. Since $\omega_j \in \cH_G$, we have $\deg \omega_j \leq \deg \Delta_G$.
    
    If $\deg \omega_j = \deg \Delta_G$, then $\omega_j = c_j \Delta_G$ for some $c_j \in K$. Then
      \[ \dif^\dagger \omega = \sum_{j=1}^n x_j \partial_{\theta_j} \omega = \sum_{j=1}^n c_j x_j \Delta_G \in \cSH_G^0 = \cH_G. \]
    Hence $\sum_{j=1}^n c_j x_j = 0$, so $c_j=0$ and $\omega = 0$.
    
    If $\deg \omega_j = \deg \Delta_G - 1$, then $\omega_j = \partial_{\ell_j} \Delta_G$ for some $\ell_j$ of degree $1$. Hence
      \[ 0 = \partial_{\dif f_i} \omega = \sum_{j=1}^n \partial_{\partial_{x_j} f_i} \omega_j
          = \sum_{j=1}^n \partial_{\ell_j \partial_{x_j} f_i} \Delta_G. \]
    By \eqref{eq:Steinberg_ann}, $\sum_{j=1}^n \ell_j \partial_{x_j} f_i \in \cI_G$. From \Cref{tab:explicit}, we may use $f_i = \sum_{j=1}^n x_j^{mi}$ for $1 \leq i \leq n-1$ and $f_n = (x_1 \cdots x_n)^{m/p}$. The $i=n$ condition becomes
    \begin{equation}\label{eq:k1_top.1}
      \sum_{j=1}^n \ell_j x_j^{-1} (x_1 \cdots x_n)^{m/p} \in \langle (x_1 \cdots x_n)^{m/p}, \sum_{j=1}^n x_j^{mi}\text{ for }1 \leq i \leq n-1\rangle.
    \end{equation}
    If $n=1$, we have $\ell_1 = c x_1$ for $c \in K$, so $\omega = c \cdot \dif \Delta_G$. We may thus assume $n \geq 2$. 

    First suppose $m/p+1 < m$. In this case, the $x_k$-degree of monomials on the left-hand side of \eqref{eq:k1_top.1} is at most $m/p+1$, while the $x_k$-degree of the generators $\sum_{j=1}^n x_j^{im}$ are at least $m$, so they cannot contribute. Hence $\sum_{j=1}^n \ell_j x_j^{-1}$ is constant. We leave it to the reader to check that this implies $\ell_j = c_j x_j$ for $c_j \in K$. The condition for $i=1$ is therefore
    \begin{equation}\label{eq:k1_top.2}
      \sum_{j=1}^n c_j x_j^m \in \langle (x_1 \cdots x_n)^{m/p}, \sum_{j=1}^n x_j^{mi}\text{ for }1 \leq i \leq n-1\rangle
    \end{equation}
    Hence we have $c \in K$ and $f \in K[x_1, \ldots, x_n]$ such that
      \[ \sum_{j=1}^n (c_j - c) x_j^m = f \cdot (x_1 \cdots x_n)^{m/p}. \]
    Since in the monomial expansion of the right-hand side of this equation every monomial is of the form $x^\beta$ for $\beta_i > 0$ for all $1 \leq i \leq n$ and $n \geq 2$, both sides of this equation must be zero. Thus $f=0$ and $c_j = c$ for all $j$. Hence $\omega = c \cdot \dif \Delta_G$.
    
    Finally, suppose $m/p + 1 \geq m$. If $p=1$, the result follows from \Cref{thm:B}, so take $p>1$. If $m>2$, then $m/p+1 \leq m/2 + 1 < m/2+m/2 = m$, contrary to our assumption, so $m=2=p$ and $m/p+1 = m$. The higher generators $\sum_{j=1}^n x_j^{mi}$ for $i \geq 2$ here cannot contribute to \eqref{eq:k1_top.1}, so we have $c \in K$ and $g \in K[x_1, \ldots, x_n]$ such that
      \[ \sum_{j=1}^n \ell_j x_j^{-1} x_1 \cdots x_n = c x_1 \cdots x_n + g (x_1^2 + \cdots x_n^2). \]
    The $x_k$-degree of the left-hand side is at most $2$, forcing the $x_k$-degree of $g$ to be $0$ for all $k$. Hence $g$ is constant. If $n \geq 3$, we see $g=0$, and the previous arguments apply to give $\omega = c \cdot \dif \Delta_G$.
    
    We are left with the case $m=p=n=2$. Here we have
      \[ \cJ_{G(2, 2, 2)} = \langle x_1^2 + x_2^2, x_1 x_2, x_1\theta_1 + x_2\theta_2, x_2\theta_1 + x_1\theta_2\rangle. \]
    It is straightforward to verify directly that the claimed elements are the only top-degree elements of $\cSH_G^1$.
  \end{proof}
\end{lemma}

\begin{example}\label{ex:G222}
  For $G = G(2, 2, 2)$, the elements of maximal total degree are spanned by $\Delta_G = x_1^2 - x_2^2$, $\dif \Delta_G/2 = x_1 \theta_1 - x_2 \theta_2$, $-\dif_2 \Delta_G/2 = x_2 \theta_1 - x_1 \theta_2$, and $\dif \dif_2 \Delta_G/4 = \theta_1\theta_2$.
\end{example}

We now give an analogue of the key argument we used in \Cref{sec:bidegree}.

\begin{lemma}\label{lem:Brunault.p}
  Let $G = G(m, p, n)$ where $p>1$. If $j \geq 1$ and $f \in K[x_j, \ldots, x_n]$, then
  \begin{align*}
    \deg f &> (n-j+1)(m(j-1)-1) \\
      &\Rightarrow\quad
        f \cdot (x_j \cdots x_n)^{m/p} \theta_j \cdots \theta_n \in \cJ_{G(m, p, n)}.
  \end{align*}
  
  \begin{proof}
    If $j=1$, then $(x_1 \cdots x_n)^{m/p} \in \cI_G \subset \cJ_G$ and the result follows. Now take $1 < j \leq n$.
    
    Suppose we have a common zero of the set of polynomials
      \[ P = \{x_i \partial_{x_i} h_{j-1}(x_j^m, \ldots, x_n^m) : i=j, \ldots, n\}. \]
    We have
      \[ \left(\sum_{i=j}^n x_i \partial_{x_i}\right) h_{j-1}(x_j^m, \ldots, x_n^m)
          = m(j-1) h_{j-1}(x_j^m, \ldots, x_n^m) = 0, \]
    so $h_{j-1}(x_j^m, \ldots, x_n^m) = 0$. Hence
    \begin{align*}
      \partial_{x_i} h_j(x_j^m, \ldots, x_n^m)
        &= \partial_{x_i} (h_j(x_j^m, \ldots, \widehat{x_i}, \ldots, x_n^m) + x_i^m h_{j-1}(x_j^m, \ldots, x_n^m)) \\
        &= m x_i^{m-1} h_{j-1}(x_j^m, \ldots, x_n^m) + x_i^m \partial_{x_i} h_{j-1}(x_j^m, \ldots, x_n^m) \\
        &= 0.
    \end{align*}
    Thus the only common zero of $P$ is $0$ by \Cref{lem:Brunault}, so $P$ forms a regular sequence in $K[x_j, \ldots, x_n]$. The Hilbert series of the quotient is hence $[m(j-1)]_q^{n-j+1}$, which has degree $(n-j+1)(m(j-1)-1)$. Thus $f \in \langle P\rangle$.
    
    By \Cref{thm:grobner_mpn}, $h_{j-1}(x_j^m, \ldots, x_n^m) (x_j \cdots x_n)^{m/p} \in \cI_G$. Also,
    \begin{align*}
      x_i \partial_{x_i} h_{j-1}&(x_j^m, \ldots, x_n^m) (x_j \cdots x_n)^{m/p} \\
        &= x_i \left(\partial_{x_i} h_{j-1}(x_j^m, \ldots, x_n^m)\right) (x_j \cdots x_n)^{m/p} \\
        &+ \frac{m}{p} h_{j-1}(x_j^m, \ldots, x_n^m) (x_j \cdots x_n)^{m/p}.
    \end{align*}
    Hence
    \begin{align*}
      x_i \dif h_{j-1}(x_j^m, \ldots, x_n^m) &(x_j \cdots x_n)^{m/p} \theta_j \cdots \widehat{\theta_i} \cdots \theta_n \\
        &= \pm \left(x_i \partial_{x_i} h_{j-1}(x_j^m, \ldots, x_n^m)\right) (x_j \cdots x_n)^{m/p} \theta_j \cdots \theta_n \\
        &\pm \frac{m}{p} h_{j-1}(x_j^m, \ldots, x_n^m) (x_j \cdots x_n)^{m/p} \theta_j \cdots \theta_n \\
        &\in \cJ_G.
    \end{align*}
    Thus $\langle P\rangle (x_j \cdots x_n)^{m/p} \theta_j \cdots \theta_n \subset \cJ_G$, and the result follows.
  \end{proof}
\end{lemma}

We may now restate and prove \Cref{thm:C}.

\begin{reptheorem}{thm:C}
  Let $G = G(m, p, n)$ with $p \neq m$ or $p=1$. Then
    \[ \bigoplus_{i+k=\ell} \cSH_{G(m, p, n)}^{i, k} \neq 0 \quad\Leftrightarrow\quad
        0 \leq \ell \leq m \binom{n}{2} + n \left(\frac{m}{p} - 1\right). \]
  Moreover, if $\ell = m \binom{n}{2} + n \left(\frac{m}{p} - 1\right)$, then
    \[ \bigoplus_{i+k=\ell} \cSH_{G(m, p, n)}^{i, k} = 
    (K + K\dif) \prod_{1 \leq i < j \leq n} (x_j^m - x_i^m) \cdot (x_1 \cdots x_n)^{m/p - 1}. \]
  
  \begin{proof}
    If $p=1$, the result follows from \Cref{thm:B} and \Cref{lem:k1_top}, so take $p>1$.
    
     Consider a monomial $x^\alpha \theta_I \not\in \cJ_G$ with $|\alpha| = i$ and $|I| = k$. Applying an element of $\fS_n$, we may take $\theta_I = \theta_{n-k+1} \cdots \theta_n$. Since $\cI_G \subset \cJ_G$, we may further take $x^\alpha \in \cA(m, p, n)$.
    
    Suppose $x^\alpha$ is obtained by $p$-contracting rows $j-1, j, \ldots, n$. Then $\alpha_{j-1} \leq \frac{m}{p} - 1$ and $\alpha_j, \ldots, \alpha_n \geq \frac{m}{p}$, so $(x_j \cdots x_n)^{m/p} \mid x^\alpha$.
    
    We first consider the case when $j \leq n-k+1$. We may write
      \[ x^\alpha = x_1^{\alpha_1} \cdots x_{n-k}^{\alpha_{n-k}} f(x_{n-k+1}, \ldots, x_n) (x_{n-k+1} \cdots x_n)^{m/p}. \]
    If $\deg f > k (m(n-k) - 1)$ then $x^\alpha \theta_I \in \cJ_G$ by \Cref{lem:Brunault.p}, so $\deg f \leq k (m(n-k) - 1)$. Since $p$-contraction starts at $j-1 \leq n-k$, $\alpha_1 + \cdots + \alpha_{n-k}$ is maximized by $p$-contracting the maximal-degree element of $\cA(m, 1, n)$, which results in $x^\beta$ with
      \[ \beta_1 = \frac{m}{p}-1, \beta_2 = \frac{m}{p}-1 + m, \ldots, \beta_{n-k} = \frac{m}{p}-1 + (n-k-1) m. \]
    Hence
      \[ \alpha_1 + \cdots + \alpha_{n-k} \leq \sum_{\ell=1}^{n-k} \left(\left(\frac{m}{p}-1\right) + (\ell-1)m\right). \]
    Combining these bounds, we find that
    \begin{align*}
      \deg \Delta_G - \deg x^\alpha \theta_I
        &\geq m \binom{n}{2} + n\left(\frac{m}{p} - 1\right) \\
        &- \sum_{\ell=1}^{n-k} \left(\left(\frac{m}{p}-1\right) + (\ell-1)m\right) \\
        &- k (m(n-k)-1) - k \frac{m}{p} - k \\
        &= m \binom{k}{2} - k.
    \end{align*}
    This difference is positive if $k \geq 2$ since $m > p > 1$, in which case $\deg \Delta_G > i+k$ as required. If $k=0$, the result follows from \Cref{thm:Steinberg}, and if $k=1$, the result follows from \Cref{lem:k1_top}.
    
    Now consider the case when $j > n-k+1$. Since $(x_j \cdots x_n)^{m/p} \mid x^\alpha$, we may write
      \[ x^\alpha = x_1^{\alpha_1} \cdots x_{j-1}^{\alpha_{j-1}} f(x_j, \ldots, x_n) (x_j \cdots x_n)^{m/p}. \]
    If $\deg f > (n-j+1)(m(j-1)-1)$ then \Cref{lem:Brunault.p} implies that
      \[ f \cdot (x_j \cdots x_n)^{m/p} \theta_j \cdots \theta_n \in \cJ_G, \]
    so $x^\alpha \theta_I \in \cJ_G$. Hence $\deg f \leq (n-j+1)(m(j-1)-1)$. Since $p$-contraction begins at row $j-1$, we have
      \[ \alpha_1 + \cdots + \alpha_{j-1} \leq \sum_{\ell=1}^{j-2} \left((m-1) + (\ell-1)m\right) + \frac{m}{p} - 1. \]
    In all, we find
    \begin{align}
      \deg \Delta_G - \deg x^\alpha \theta_I \nonumber
        &\geq m \binom{n}{2} + n\left(\frac{m}{p} - 1\right) \nonumber \\
        &- \sum_{\ell=1}^{j-2} \left((m-1) + (\ell-1)m\right) - \frac{m}{p} + 1 \nonumber \\
        &- (n-j+1)(m(j-1)-1) - (n-j+1) \frac{m}{p} - k \nonumber \\
        &= \frac{m}{p} (j-2) + m \binom{n-j+1}{2} - k \label{eq:thm:C.1}.
    \end{align}
    
    Suppose $n-j \geq 2$. Then since $p>1$, the bound in \eqref{eq:thm:C.1} becomes
    \begin{align*}
      \frac{m}{p} (j-2) &+ \frac{m}{2} (n-j+1)(n-j) - k \\
        &\geq \frac{m}{p} \left((j-2) + (n-j+1)(n-j)\right) - k \\
        &\geq \frac{m}{p} \left((j-2) + 2(n-j+1)\right) - k \\
        &> \frac{m}{p} \left(j-2+n-j+2\right) - n \\
        &= n \left(\frac{m}{p} - 1\right) \geq 0.
    \end{align*}
    Thus, as before we have $\deg \Delta_G > i+k$.
    
    We are left with the cases $j \in \{n-1, n, n+1\}$. When $n=1$ or $n=2$, the result follows from \Cref{lem:k1_top} and \Cref{cor:dif_difdagger}, so take $n \geq 3$. If $j=n+1$, then since $m > p$, the bound in \eqref{eq:thm:C.1} is
      \[ \frac{m}{p} (n-1) - k \geq 2(n-1)-k \geq 2(n-1)-n = n-2 \geq 1. \]
    If $j=n-1$, the bound becomes
    \begin{align*}
      \frac{m}{p} (n-3) + m - k
        &\geq 2(n-3) + 4 - n = n-2 \geq 1.
    \end{align*}
    
    Finally, if $j=n$, we have
    \begin{align*}
      \frac{m}{p} (n-2) - k
        &\geq 2(n-2) - k \geq 2(n-2)-n = n-4,
    \end{align*}
    so the result follows when $n \geq 5$. We are left with the cases $n \in \{3, 4\}$. If $n=3$, we have $2(n-2) - k = 2 - k$, so the result holds at $k=0, 1$. It also holds at $k=3$ by \Cref{cor:top_xdeg}, so it holds at $k=2$ by \Cref{cor:dif_difdagger}. If $n=4$, we have $2(n-2) - k = 4 - k$, so the result holds at $k=0, 1, 2, 3$. It also holds at $k=4$ by \Cref{cor:top_xdeg}. This completes the proof.
  \end{proof}
\end{reptheorem}

\begin{remark}
  The preceding arguments handle the vast majority of the groups with $m=p$ as well. One finds that the result holds except possibly when $k \in \{n-2, n-1\}$ and when $p$-contraction starts at the $(n-1)$st row. Note that the ``extra'' elements $\dif_2 \Delta_{G(2, 2, 2)}$ and $\dif \dif_2 \Delta_{G(2, 2, 2)}$ from \Cref{ex:G222} are of this form, though they are $\det$-isotypic.
\end{remark}

\begin{center}
\begin{sidewaystable}
  \centering
  \begin{tabular}{c|c|c|c|c|c|c|c}
    $G$
      & $f_i$
      & $\Delta_G$
      & $\deg \Delta_G$
      & $\dif_i$
      & $e_i^*$
      & $\Delta_G^*$
      & $\deg \Delta_G^*$ \\
      
    \toprule
    \begin{textcolumn}
    ${\scriptstyle G(1, 1, n)}$ \\
    ${\scriptstyle = \fS_n}$
    \end{textcolumn}
      & ${\scriptstyle \prod_{j=1}^n x_j^i}$
      & ${\scriptstyle \prod\limits_{1 \leq i < j \leq n} (x_j - x_i)}$
      & ${\scriptstyle \binom{n}{2}}$
      & ${\scriptstyle \sum\limits_{j=1}^n \partial_{x_j^i}} \theta_j$
      & ${\scriptstyle i}$
      & ${\scriptstyle \prod\limits_{1 \leq i < j \leq n} (x_j - x_i)}$
      & ${\scriptstyle \binom{n}{2}}$ \\
    
    \midrule
    \begin{textcolumn}
    ${\scriptstyle G(2, 1, n)}$ \\
    ${\scriptstyle = \fB_n}$
    \end{textcolumn}
      & ${\scriptstyle \sum\limits_{j=1}^n x_j^{2i}}$
      & \begin{textcolumn}
        ${\scriptstyle \prod\limits_{1 \leq i < j \leq n} (x_j^2 - x_i^2)}$ \\
        ${\scriptstyle \cdot x_1 \cdots x_n}$
        \end{textcolumn}
      & ${\scriptstyle n^2}$
      & ${\scriptstyle \sum\limits_{j=1}^n \partial_{x_j^{2i-1}} \theta_j}$
      & ${\scriptstyle 2i-1}$
      & \begin{textcolumn}
         ${\scriptstyle \prod\limits_{1 \leq i < j \leq n} (x_j^2 - x_i^2)}$ \\
         ${\scriptstyle \cdot x_1 \cdots x_n}$
         \end{textcolumn}
      & ${\scriptstyle n^2}$ \\
    
    \midrule
    \begin{textcolumn}
    ${\scriptstyle G(2, 2, n)}$ \\
    ${\scriptstyle = \fD_n}$
    \end{textcolumn}
      & \begin{textcolumn}
          ${\scriptstyle \sum\limits_{j=1}^n x_j^{2i}}$ \\
          \scriptsize{(for ${\scriptstyle 1 \leq i \leq n-1}$)} \\
          \  \\
          ${\scriptstyle x_1 \cdots x_n}$ \\
          \scriptsize{(for ${\scriptstyle i=n}$)}
          \end{textcolumn}
      & ${\scriptstyle \prod\limits_{1 \leq i < j \leq n} (x_j^2 - x_i^2)}$
      & ${\scriptstyle n(n-1)}$
      & \begin{textcolumn}
          ${\scriptstyle \sum\limits_{j=1}^n \partial_{x_j^{2i-1}} \theta_j}$ \\
          \scriptsize{(for ${\scriptstyle 1 \leq i \leq n-1}$)} \\
          \  \\
          ${\scriptstyle \sum\limits_{j=1}^n \partial_{x_1 \cdots \widehat{x}_j \cdots x_n} \theta_j}$ \\
          \scriptsize{(for ${\scriptstyle i=n}$)}
          \end{textcolumn}
      & \begin{textcolumn}
          ${\scriptstyle 2i-1}$ \\
          \scriptsize{(for ${\scriptstyle 1 \leq i \leq n-1}$)} \\
          \  \\
          ${\scriptstyle n-1}$ \\
          \scriptsize{(for ${\scriptstyle i=n}$)}
          \end{textcolumn}
      & ${\scriptstyle \prod\limits_{1 \leq i < j \leq n} (x_j^2 - x_i^2)}$
      & ${\scriptstyle n(n-1)}$ \\

    \midrule
    \begin{textcolumn}
    ${\scriptstyle G(m, m, 2)}$ \\
    ${\scriptstyle = \fDih_{2m}}$ \\
    \scriptsize{for ${\scriptstyle m \geq 2}$}
    \end{textcolumn}
      & \begin{textcolumn}
          ${\scriptstyle x_1^m + x_2^m,}$ \\
          ${\scriptstyle x_1 x_2}$
          \end{textcolumn}
      & ${\scriptstyle x_2^m - x_1^m}$
      & ${\scriptstyle m}$
      & \begin{textcolumn}
          ${\scriptstyle \partial_{x_1} \theta_1 + \partial_{x_2} \theta_2,}$ \\
          ${\scriptstyle \partial_{x_2^{m-1}} \theta_1 + \partial_{x_1^{m-1}} \theta_2}$
          \end{textcolumn}
      & \begin{textcolumn}
          ${\scriptstyle 1,}$ \\
          ${\scriptstyle m-1}$
          \end{textcolumn}
      & ${\scriptstyle x_2^m - x_1^m}$
      & ${\scriptstyle m}$ \\

    \midrule
    \begin{textcolumn}
      ${\scriptstyle G(m, 1, n)}$ \\
      \scriptsize{for ${\scriptstyle m>1}$}
    \end{textcolumn}
      & ${\scriptstyle \sum\limits_{j=1}^n x_j^{mi}}$
      & \begin{textcolumn}
        ${\scriptstyle \prod\limits_{1 \leq i < j \leq n} (x_j^m - x_i^m)}$ \\
        ${\scriptstyle \cdot (x_1 \cdots x_n)^{m-1}}$
        \end{textcolumn}
      & ${\scriptstyle m \binom{n}{2} + n(m-1)}$
      & ${\scriptstyle \sum\limits_{j=1}^n \partial_{x_j^{(i-1)m+1}} \theta_j}$
      & ${\scriptstyle (i-1)m+1}$
      & \begin{textcolumn}
         ${\scriptstyle \prod\limits_{1 \leq i < j \leq n} (x_j^m - x_i^m)}$ \\
         ${\scriptstyle \cdot x_1 \cdots x_n}$
         \end{textcolumn}
      & ${\scriptstyle m \binom{n}{2} + n}$ \\

    \midrule
    \begin{textcolumn}
    ${\scriptstyle G(m, p, n)}$ \\
    \scriptsize{for ${\scriptstyle p \neq m}$}
    \end{textcolumn}
      & \begin{textcolumn}
          ${\scriptstyle \sum\limits_{j=1}^n x_j^{mi}}$ \\
          \scriptsize{(for ${\scriptstyle 1 \leq i \leq n-1}$)} \\
          \ \\
          ${\scriptstyle (x_1 \cdots x_n)^{m/p}}$ \\
          \scriptsize{(for ${\scriptstyle i=n}$)}
          \end{textcolumn}
      & \begin{textcolumn}
          ${\scriptstyle \prod\limits_{1 \leq i < j \leq n} (x_j^m - x_i^m)}$ \\
          ${\scriptstyle \cdot (x_1 \cdots x_n)^{m/p - 1}}$
          \end{textcolumn}
      & ${\scriptstyle m \binom{n}{2} + n\left(\frac{m}{p} - 1\right)}$
      & ${\scriptstyle \sum\limits_{j=1}^n \partial_{x_j^{(i-1)m+1}} \theta_j}$
      & ${\scriptstyle (i-1)m + 1}$
      & \begin{textcolumn}
          ${\scriptstyle \prod\limits_{1 \leq i < j \leq n} (x_j^m - x_i^m)}$ \\
          ${\scriptstyle \cdot x_1 \cdots x_n}$
          \end{textcolumn}
      & ${\scriptstyle m \binom{n}{2} + n}$ \\
    
    \midrule
    ${\scriptstyle G(m, m, n)}$
      & \begin{textcolumn}
          ${\scriptstyle \sum\limits_{j=1}^n x_j^{mi}}$ \\
          \scriptsize{(for ${\scriptstyle 1 \leq i \leq n-1}$)} \\
          \  \\
          ${\scriptstyle x_1 \cdots x_n}$ \\
          \scriptsize{(for ${\scriptstyle i=n}$)}
          \end{textcolumn}
      & ${\scriptstyle \prod\limits_{1 \leq i < j \leq n} (x_j^m - x_i^m)}$
      & ${\scriptstyle m \binom{n}{2}}$
      & \begin{textcolumn}
          ${\scriptstyle \sum\limits_{j=1}^n \partial_{x_j^{(i-1)m + 1}} \theta_j}$ \\
          \scriptsize{(for ${\scriptstyle 1 \leq i \leq n-1}$)} \\
          \  \\
          ${\scriptstyle \sum\limits_{j=1}^n \partial_{(x_1 \cdots \widehat{x}_j \cdots x_n)^{m-1}} \theta_j}$ \\
          \scriptsize{(for ${\scriptstyle i=n}$)}
          \end{textcolumn}
      & \begin{textcolumn}
          ${\scriptstyle (i-1)m + 1}$ \\
          \scriptsize{(for ${\scriptstyle 1 \leq i \leq n-1}$)} \\
          \  \\
          ${\scriptstyle (n-1)(m-1)}$ \\
          \scriptsize{(for ${\scriptstyle i=n}$)}
          \end{textcolumn}
      & ${\scriptstyle \prod\limits_{1 \leq i < j \leq n} (x_j^m - x_i^m)}$
      & ${\scriptstyle m \binom{n}{2}}$
  \end{tabular}
  \vspace{0.5em}
  \caption{Summary of explicit descriptions of the basic invariants $f_i$, the Vandermondians $\Delta_G$ and their degrees, the generalized exterior derivatives $\dif_i$, the co-exponents $e_i^*$, and the co-Vandermondians $\Delta_G^*$ and their degrees for all non-exceptional pseudo-reflection groups $G$. Here $f_i$ ranges over $1 \leq i \leq n$, and $d_i$ ranges over $1 \leq i \leq r$ where $r=n$ except $r=n-1$ when $G = \fS_n$.}
  \label{tab:explicit}
\end{sidewaystable}
\end{center}

\begin{center}
\begin{sidewaystable}
  \centering
  \begin{tabular}{c|c|c|c}
    $G$ & $\Hilb(\cSH_G; 1, z)$ & $\Hilb(K[\partial_{x_1}, \ldots, \partial_{x_n}]\cSH_G^{\det}; 1, z)$ \\
    
    \toprule
    \toprule
    $\fS_3$ & $z^2 + 6z + 6$ & (same) \\
    
    \midrule
    $\fS_4$ & $z^3 + 14z^2 + 36z + 24$ & (same) \\
    
    \midrule
    $\fS_5$ & $z^4 + 30z^3 + 150z^2 + 240z + 120$ & (same) \\
    
    \midrule
    $\fS_6$ & $z^5 + 62z^4 + 540z^3 + 1560z^2 + 1800z + 720$ & (same) \\
    
    \midrule
    $\fB_4$ & $z^4 + 80z^3 + 464z^2 + 768z + 384$ & (same) \\
    
    \midrule
    $\fB_5$ & $z^5 + 242z^4 + 2640z^3 + 8160z^2 + 9600z + 3840$ & (same) \\
    
    \midrule
    $\fD_2$ & $z^2 + 4z + 4$ & (same) \\

    \midrule
    $\fD_3$ & $z^3 + 14z^2 + 36z + 24$ & (same) \\
    
    \midrule
    $\fD_4$ & $z^4 + 48z^3 + 240z^2 + 384z + 192$ & $z^4 + 46z^3 + 238z^2 + 384z + 192$ \\
    
    \midrule
    $\fD_5$ & $z^5 + 162z^4 + 1440z^3 + 4160z^2 + 4800z + 1920$ & $z^5 + 147z^4 + 1405z^3 + 4140z^2 + 4800z + 1920 $ \\
    
    \midrule
    $H_3$ & $z^3 + 62z^2 + 180z + 120$ & (same) \\
    
    \midrule
    $F_4$ & $z^4 + 244z^3 + 1396z^2 + 2304z + 1152$ & $z^4 + 220z^3 + 1372z^2 + 2304z + 1152$ \\
    
    \midrule
    $G(3, 1, 2)$ & $4z^2 + 21z + 18$ & (same) \\

    \midrule
    $G(4, 1, 2)$ & $9z^2 + 40z + 32$ & (same) \\
    
    \midrule
    $G(5, 1, 2)$ & $16z^2 + 65z + 50$ & (same) \\
    
    \midrule
    $G(5, 1, 3)$ & $64z^3 + 665z^2 + 1350z + 750$ & (same) \\
    
    \midrule
    $G(3, 1, 4)$ & $16z^4 + 609z^3 + 2862z^2 + 4212z + 1944$ & (same) \\

    \midrule
    $G(4, 1, 4)$ & $81z^4 + 2320z^3 + 9920z^2 + 13824z + 6144$ & (same) \\
 
    \midrule
    $G(4, 2, 4)$ & $76z^4 + 1451z^3 + 5408z^2 + 7104z + 3072$ & $z^4 + 544z^3 + 3616z^2 + 6144z + 3072$ \\
 
   \midrule
   $G(4, 4, 4)$ & $33z^4 + 416z^3 + 1920z^2 + 3072z + 1536$ & $z^4 + 286z^3 + 1822z^2 + 3072z + 1536$
  \end{tabular}
  \vspace{0.5em}
  \caption{Computer calculations of Hilbert series of the terms in \eqref{eq:thm:A}. The $q=1$ specialization has been taken to make the output manageable. \eqref{eq:thm:A} holds for a given $G$ if and only if the two columns agree. In each case, \eqref{eq:lem:B} and \eqref{eq:lem:C.1} indeed hold, except \eqref{eq:lem:B} fails for $G(4, 2, 4)$ and $G(4, 4, 4)$. Calculations were done using SageMath \cite{sagemath}. The exceptional group calculations were done with \cite{MR926338}.}
  \label{tab:calcs}
\end{sidewaystable}
\end{center}

\clearpage

%\section{Acknowledgements}
%
%(To be written.)

\bibliography{refs}{}
\bibliographystyle{acm}

\end{document}